\documentclass[final,3p,times]{elsarticle}

\usepackage[normalem]{ulem}

\usepackage{amssymb}
\usepackage{amsmath}
\usepackage{amsthm}
\usepackage{subeqnarray}
\usepackage{mathrsfs}
\usepackage{rotating}
\usepackage{enumitem}
\setlist[enumerate]{noitemsep, topsep=2pt}
\setlist[itemize]{noitemsep, topsep=2pt}
\usepackage{bm}             
\usepackage{setspace}
\usepackage{numprint}

\usepackage{algorithmic}
\usepackage[linesnumbered,ruled]{algorithm2e}

\usepackage[breaklinks,bookmarks=false]{hyperref}
\hypersetup{colorlinks, linkcolor=blue, citecolor=blue,
  urlcolor=blue, pagecolor=blue}

\usepackage{cleveref}
\usepackage{amstext}

\biboptions{sort&compress}

\newcommand{\si}[1]{\,\mathrm{#1}}
\newcommand{\np}[1]{\,\numprint{#1}}

\newtheorem{theorem}{Theorem}[section]

\theoremstyle{remark}
\newtheorem{remark}[theorem]{Remark}
\newtheorem{definition}[theorem]{Definition}

\newcommand\bifont {\bm}

\newcommand\krt{{\mathscr T}}

\newcommand\bkn{\bifont{n}}
\newcommand\bku{\bifont{u}}
\newcommand\bkv{\bifont{v}}

\newcommand{\vertical}[1]{\begin{sideways}{#1}\end{sideways}}

\newcommand\pd {\partial}
\newcommand{\pdt}[1]{ \pd_t{#1} }
\newcommand{\pdd}[2]{ \pd_{#2}{#1} }

\newcommand\R  {\mathbb{R}}

\newcommand\Om {{\Omega}}
\newcommand\oO {\overline{\om}}

\newcommand\gom{{\Gamma}}

\newcommand\dd {\mathrm{d}}
\newcommand\dt {{\,\dd t}}

\newcommand\dx {{\,\dd x}}

\newcommand\pp {{\bifont{p}}}
\newcommand\qq {{q}}

\newcommand\Nhm {{N_{h,m}}}
\newcommand\Mhm {{M_{h,m}}}

\newcommand\bbw {\bar{{\w}}}
\newcommand\wht {{\w}_{h\tau}}

\newcommand\qht {{q}_{h\tau}}
\newcommand\Wqht {\widehat{q}_{h\tau}}

\newcommand\wm {{\w}_{h\tau}^m}
\newcommand\wmk {{\w}_{h\tau}^{m,k}}
\newcommand\wmkM {{\w}_{h\tau}^{m,k-1}}
\newcommand\wmkMs {{\w}_{h\tau}^{m,k-s}}
\newcommand\wmN {{\w}_{h\tau}^{m,0}}
\newcommand\wmNN {{\w}_{h\tau}^{m-1}}
\newcommand\buk {{\bku}_{h\tau}^{k}}

\newcommand\vv {\bifont{v}}

\newcommand\TT {\mathrm{T}}
\newcommand\Th {\krt_h}
\newcommand\Thm{{\krt_{h,m}}}
\newcommand\ThmM{{\krt_{h,m-1}}}
\newcommand\ThmP{{\krt_{h,m+1}}}
\newcommand\Thn{{\krt_{h,0}}}

\newcommand\Fhm{{{\mathcal F}_{h,m}}}
\newcommand\FhmI{{{\mathcal F}_{h,m}^I}}
\newcommand\FhmB{{{\mathcal F}_{h,m}^B}}

\newcommand{\norm}[2]{\left\|#1\right\|_{#2} }
\newcommand{\normP}[3]{\left\|#1\right\|_{#2}^{#3} }

\newcommand\ahm {{a_{h,m}}}
\newcommand\ahmL {{a_{h,m}^{\scriptscriptstyle\mathrm{L}}}}
\newcommand\Ahm {{A_{h,m}}}
\newcommand\AhmL {{A_{h,m}^{\scriptscriptstyle\mathrm{L}}}}

\newcommand\ahmt {{\tilde{a}_{h,m}}}

\newcommand\om{{\Omega}}
\newcommand\cJ{\mathcal{J}}

\newcommand\plus{^{\scriptscriptstyle (+)}}
\newcommand\minus{^{\scriptscriptstyle (-)}}

\newcommand\Plus{^{\scriptscriptstyle+}}
\newcommand\Minus{^{\scriptscriptstyle-}}
\newcommand\PMinus{^{\scriptscriptstyle\pm}}

\newcommand\bF{\bifont{f}}
\newcommand\bFv{\vec{\ff}}

\newcommand\bA{\bifont{A}}
\newcommand\bB{\bifont{B}}
\newcommand\bAv{\vec{\bA}}
\newcommand\bP {\bifont{P}}
\newcommand\PP {\bP}
\newcommand\bR {\bifont{R}}
\newcommand\bRv {\vec{\bR}}
\newcommand\mK {\bifont{K}}
\newcommand\mKv {\widehat{\mK}}

\renewcommand\bkn{\vec{n}}
\newcommand\bknK {{\bkn_{\!\!\;\scriptscriptstyle K}}}

\newcommand\etaAm {\eta^{\scriptscriptstyle\mathrm{A}}_{m}}
\newcommand\etaSm {\eta^{\scriptscriptstyle\mathrm{S}}_{m}}
\newcommand\etaTm {\eta^{\scriptscriptstyle\mathrm{T}}_{m}}

\newcommand\TOL {\omega}

\newcommand\jump[1]{[\![{#1}]\!]}
\newcommand{\tjump}[1]{\left\{{#1}\right\}}
\newcommand\aver[1]{\left\langle{#1}\right\rangle}
\newcommand{\Lsp}[2]{\left({#1},{#2}\right)_{\Omega}}
\newcommand{\LSP}[3]{\left({#1},{#2}\right)_{#3}}

\newcommand\sumk{{ \sum_{K\in\Thm} }}

\newcommand\sumK[2]{{ \sum_{K\in\Thm}\left( {#1}, {#2} \right)_{K} }}
\newcommand\sumKK[3]{{ \sum_{K\in\Th}\left( {#1}, {#2} \right)_{#3} }}
\newcommand\sumF[2]{{ \sum_{\gamma\in\Fhm}\left( {#1}, {#2} \right)_{\gamma} }}

\newcommand\sumFB[2]{{ \sum_{\gamma\in\FhmB}\left( {#1}, {#2} \right)_{\gamma} }}

\newcommand\DoF{\mathrm{DoF}}

\newcommand\w {{\bifont{w}}}
\newcommand\ff {{\bifont{f}}}
\newcommand\bg {{\bifont{g}}}
\newcommand\bk {{\bifont{k}}}

\newcommand\dK {\partial K}

\newcommand\press {\mathrm{p}}
\newcommand\pressE {\mathrm{P}}
\newcommand\ener {e}
\newcommand\ccP{c_{\press}}
\newcommand\ccV{c_{\mathrm{v}}}
\newcommand\stress{\tau}

\newcommand\temp{\theta}
\newcommand\tempP{\Theta}

\newcommand\pK {p_K}

\newcommand\mbS{\bifont{S}}

\newcommand\Smhp {S_{m,h,\pp}}
\newcommand\SmhpP {S_{m,h,\pp+1}}
\newcommand\bSmhp {\mbS_{m,h,\pp}}
\newcommand\bSmhpP {\mbS_{m,h,\pp+1}}
\newcommand\bSmhpM {\mbS_{m-1,h,\pp}}
\newcommand{\vas} {\varphi}
\newcommand{\bvas} {\pmb{\varphi}}
\newcommand{\bpsi} {\pmb{\psi}}

\newcommand\bShpq {\mbS_{h,\pp}^{\tau,q}}
\newcommand\bShpqS {\mbS_{h,\pp+1}^{\tau,\qq}}
\newcommand\bShpqT {\mbS_{h,\pp}^{\tau,\qq+1}}

\newcommand\etam {\eta_m^{\mathrm{int}}}
\newcommand \Mir{{\mathscr M}}  
\newcommand \BC{{\mathscr B}}
\newcommand \Ext{{\mathscr E}}

\newcommand \hpAMA{{$hp$-AMA}}

\usepackage[usenames]{color}

\newcommand{\Smaz}[1]{{}}

\begin{document}


\begin{frontmatter}



\title{Non-hydrostatic mesoscale atmospheric modeling by the anisotropic mesh adaptive discontinuous Galerkin method\tnoteref{label1}}
\tnotetext[label1]{This work was supported by grant No. 20-01074S of the Czech Science Foundation.
  The author acknowledges also the membership in the Ne{\v c}as Center for Mathematical Modeling
  ncmm.karlin.mff.cuni.cz.}

\author[Prague]{V{\'\i}t Dolej{\v s}{\'\i}}
\ead{dolejsi@karlin.mff.cuni.cz}


\address[Prague]{Charles University, Faculty of Mathematics and Physics,
Sokolovsk\'a 83, 186 75 Praha, Czech Republic}


\begin{abstract}
  We deal with non-hydrostatic mesoscale atmospheric modeling using
  the fully implicit space-time discontinuous Galerkin method in combination
  with the anisotropic $hp$-mesh adaptation technique.
  The time discontinuous approximation allows the treatment of different meshes at different
  time levels in a natural way which can significantly reduce the number of degrees
  of freedom.
  The presented approach generates
  a sequence of triangular meshes consisting of possible anisotropic elements and varying
  polynomial approximation degrees such that the interpolation error is below the given tolerance
  and the number of degrees of freedom at each time step is minimal. We describe 
  the discretization of the problem together with several implementation issues
  related to the treatment of boundary conditions,
  algebraic solver and adaptive choice of the size of the time
  steps.
  The computational performance of the proposed method is demonstrated
  on several benchmark problems.
\end{abstract}
\end{frontmatter}
  

\section{Introduction}

The numerical solution of the compressible Euler and/or Navier-Stokes
equations has become the standard tool for the simulation of the atmosphere
in the last decades, see e.g. the reviews in \cite{SteppelrALL_MAP03} or
\cite{GiraldoRostelli_JCP08}.  
Due to the large scale simulation, it is advantageous to use higher order discretization
schemes since they significantly reduce the  number of degrees of freedom ($\DoF$)
needed to achieve the required accuracy. It is possible to use
the spectral element (SE) method \cite{Giraldo_JCP05,Giraldo_MWR04}
or the finite volume (FV) approach with higher-order reconstruction schemes
as ENO, WENO or TENO,
cf.~\cite{NavasMontilla_JCP23} and the references cited therein.
In the recent years, large attention has been devoted to the discontinuous Galerkin (DG) methods
which employ the advantageous properties of SE and FV methods,
cf.~\cite{RestelliGiraldo_SISC09,YelashALL_JCP14,BustoALL_CF20,WaruszewskiALL_JCP22,JiangALL_JCP22}
for example.
DG methods employ piecewise  polynomial but discontinuous approximation, they
provide high accuracy and flexibility,
local conservation, stability, and parallel efficiency.
However, the disadvantage is a large number of $\DoF$ in comparison to FV and SE methods.

Due to the high-order spatial approximation, the use of an explicit time discretization
is very inefficient, since the size of the time steps is strongly restricted by the stability
condition, cf.~\cite{Giraldo-book}.
Therefore, several variants of semi-implicit time discretization techniques have been
developed, e.g.,
the implicit-explicit (IMEX) formulation in \cite{RestelliGiraldo_SISC09, GiraldoALL_SISC10},
its recent improvement in \cite{ReddyALL_JCP23},
the IMEX evolution Galerkin method in \cite{LukacovaALL_SISC04,YelashALL_JCP14},
the positivity-preserving well-balanced central discontinuous Galerkin method in
\cite{JiangALL_JCP22} or the techniques based on the Eulerian-Lagrangian approach
in \cite{BustoALL_CF20}.

In the last two decades, we have developed, analyzed and tested the discontinuous Galerkin method
for solving the compressible Navier-Stokes equations with applications
in aerodynamics, including supersonic, transonic, subsonic, and low-Mach number flows,
cf.~\cite{DGM-book}.
The prominent scheme seems to us to be the fully implicit space-time discontinuous Galerkin (STDG)
method, which uses a discontinuous piecewise polynomial approximation with
respect to both spatial and temporal coordinates.
In our experience, the time discontinuous Galerkin method is more stable than IMEX techniques
and
it can naturally handle different meshes at different time levels. However,
these advantages are
compensated by a higher number of $\DoF$ compared to multi-step or
diagonally implicit Runge-Kutta methods.

A fully implicit time discretization leads to the need
to solve a nonlinear algebraic system 
at each time step.  We employ an iterative Newton-like method whose first step
corresponds to our earlier semi-implicit scheme from \cite{CiCP}.
Therefore, the additional iterative steps improve the
accuracy (and stability) of the computation. In practice, only a few iterative steps
are necessary to reduce the algebraic errors.

In the last decade, we have dealt with the {\em anisotropic $hp$-mesh adaptation} ({\hpAMA}) method
which constructs meshes consisting of triangles of different size, shape and orientation
together with varying polynomial approximation degrees for each mesh element.
The idea is to minimize the number of $\DoF$ such that
an estimate of the interpolation error is below the given tolerance.
This anisotropic mesh adaptation
technique is more general than the approach from \cite{WaruszewskiALL_JCP22} which allows
the movement of grid nodes but not local refinement/coarsening.
On the other hand, the adaptive technique from \cite{OrlandoALL_JCAM23} offers
the local mesh refinement/de-refinement but the initial mesh structure is fixed and then
the capturing of phenomena that are not parallel to the boundary domain is not optimal.
The theoretical as well as practical
issues of the {\hpAMA} method are summarized in \cite{AMA-book}.

The STDG discretization can naturally handle grids generated by
the {\hpAMA} method since it is
a one-step method and the approximate solutions at different time levels
are weakly connected together
by a penalty term.
In \cite{AMAtdp}, we developed the {\hpAMA} approach for solving time-dependent
problems and demonstrated its performance by solving the Burgers equation,
the isentropic vortex propagation, the Kelvin-Helmholtz instability and
single-ring infiltration in porous media flow.
The computational benefit of anisotropic meshes
for the numerical solution of unsteady flow problems
has been demonstrated in many works, e.g., 
\cite{BelmeDervieuxAlauzet_JCP12,AlauzetLoseille_JCP18,FreyAlauzet,HabashiAll} and
the references cited therein.

The goal of this paper is to apply the STDG and {\hpAMA} methods to
the numerical simulation of non-hydrostatic mesoscale atmospheric problems.
The presented algorithm controls the errors arising from an inexact solution of nonlinear
algebraic systems and chooses adaptively the size of the time steps whereas
the polynomial approximation degree with respect to time is kept fixed.
Hence, we can speak of $hp\tau$-adaptation.
Moreover, at each time step, the algorithm evaluates an estimate of the interpolation
error. If the error estimate exceeds the prescribed tolerance,
a completely new (possibly anisotropic)
mesh with varying polynomial approximation degrees with respect to the space is generated.
Attention is paid to the realization
of the boundary conditions and the choice of the physical quantity to be controlled.
The particular aspects of the presented numerical scheme are justified by numerical experiments.

The outline of the rest of the paper is as follows.
In Section~\ref{sec:problem}, we present the governing
equations with the constitutive relations and the initial/boundary conditions.
In Section~\ref{sec:DGM}, we introduce the discretization of the problem including
the realization of the boundary conditions. Furthermore, Section~\ref{sec:sol} presents
the solution strategy of the arising algebraic system with an adaptive stopping
criterion and the adaptive choice of the time step.
Section~\ref{sec:AMA} contains a description of the {\hpAMA} technique and
its application to atmospheric modeling.
The performance of the proposed method is demonstrated in Section~\ref{sec:num} by
numerical experiments for several benchmark examples from \cite{GiraldoRostelli_JCP08}.
Finally, we give several concluding remarks in Section~\ref{sec:concl}.

\section{Problem formulation}
\label{sec:problem}

Let $\Om\subset\R^2$ be the computational domain, $T>0$ the physical time  to be reached,
and $Q_T:= \Om\times(0,T)$. We use the following notation: $\rho$ is the density,
$\press$ is the pressure, 
$\vv=(v_1,v_2)^\TT$ is the velocity vector, 
$\temp$ is the absolute temperature and
$\ener = \rho\ccV\temp + \rho |\bkv|^2/2$ is the energy per unit volume
including the interior and kinetic energies (and excluding the gravitational energy).
Moreover, 
$\ccV >0$ and 
$\ccP > 0$ are the specific heat capacities at constant volume and pressure, respectively, 
$\kappa= \ccP/\ccV > 1$ is the  Poisson adiabatic constant,
$R=\ccP-\ccV $ is the gas constant,
$\mu$ is the dynamic viscosity, $g=9.81\,\mathrm{m}\cdot\mathrm{s}^{-2}$ is the gravity constant,
$\bk=(k_1,k_2)^\TT$ is the upward pointing unit vector, $\Pr$ is the Prandtl number,
and $\stress_{i j},\ i,j=1,2$ denote the components of the viscous part of the stress tensor.
The symbols $\pdt{}$ and $\pdd{}{x_i}$ denote the partial derivatives with respect to
$t$ and $x_i,\ i=1,2$, respectively.

We consider the system of the compressible Navier-Stokes equations with the
gravity forces written in the conservative form \cite{GiraldoRostelli_JCP08,OrlandoALL_JCAM23}
\begin{align} \label{eq:NS}
  \pdt{\w} + \nabla \cdot \bFv(\w)  - \nabla \cdot \bRv(\w, \nabla \w) 
  = \bg \quad\mbox{in} \ Q_T,
\end{align}
where $\nabla$ and $\nabla\cdot$ denote the gradient and divergence operators, respectively.
The unknown state vector $\w$, the convective flux $\bFv = (\bF_1, \bF_2)$,
the diffusive one $\bRv = (\bR_1, \bR_2)$
and the exterior forces $\bg$ are given by 
\begin{align}
  \label{eq:w}
  \w &=(\rho,\,\rho v_1,\rho v_{2},\,\ener)^{\TT},\qquad
  \bF_i = (\rho v_i,\,\rho v_i v_1 + \delta_{i1}\press,\,\rho v_i v_2 + \delta_{i2}\press,\,
  (\ener+\press)\,v_i)^{\TT},\qquad i=1,2, \\
   \bR_i &= \left(0,\,\stress_{i1},\,\stress_{i2},
   \sum\nolimits_{k=1}^2 \stress_{ik}v_i
   + \frac{\mu\,\ccP}{\Pr }
  \pdd{\temp}{x_i}\right)^{\TT},\ i=1,2, \qquad
  \bg = \left(0,\, -\rho g k_1,\, -\rho g k_2,\, -\rho g \bk\cdot\bkv\right)^{\TT}, \notag
\end{align}
where $\delta_{i j}$ is the Kronecker delta and the symbol $\cdot$ denotes the scalar product.
The relations \eqref{eq:NS} -- \eqref{eq:w} are accompanied by the 
constitutive relations of a perfect gas and the viscous part of the stress tensor as
\begin{align}
  \label{eq:tau}
  \press=R \rho \theta \qquad\mbox{and} \qquad
  \stress_{ij} =\mu  \left( 
 \pdd{v_i}{x_j} + \pdd{v_j}{x_i}\right)
 - 
 \tfrac23  \delta_{ij} \nabla\cdot \bkv,
 \quad i,j=1,2,
\end{align}
respectively.
In atmospheric applications, the potential temperature $\tempP$ is the preferred physical
quantity given as
\begin{align}
  \label{eq:tempP}
  \tempP = {\temp}/{\pressE},\qquad\mbox{where}\quad 
  \pressE = \left({\press}/{\press_0}\right)^{(\kappa-1)/\kappa}
\end{align}
is the Exner pressure and $\press_0 = 10^5\,\si{Pa}$ is the reference pressure.

The viscous fluxes can be written in the form
\begin{align}
  \label{Kij}
  \bRv(\w,\nabla \w) = \mKv(\w)\nabla \w \qquad \Longleftrightarrow\qquad
  \bR_i(\w,\nabla\w)=\sum\nolimits_{j=1}^{2}\mK_{i,j}(\w) \pdd{\w}{x_j},
  \quad i=1,2
\end{align}
where $\mKv$ is the ``matrix of matrices'' (fourth-order tensor) whose entries are
$4\times 4$ matrices $\mK_{i,j}(\cdot)$ dependent on $\w$ (and independent of $\nabla \w$).
Their explicit form of $\mK_{i,j}(\cdot)$, $i,j=1,2$ can be found 
e.g., in \cite{CiCP} or  \cite[Section~9.1]{DGM-book}.

The term $\bg$ representing the gravity forces can be expressed as
\begin{align}
  \label{bg}
  \bg = \bB \w,\qquad \bB =
  \begin{pmatrix}
    0 & 0 & 0 & 0 \\
    0 & 0 & 0 & 0 \\
    -g & 0 & 0 & 0 \\
    0 & 0 & -g & 0 \\
  \end{pmatrix}, 
\end{align}
where $g$ is the gravity constant (note that $\bk=(k_1,k_2)^\TT=(0,1)^\TT$).

The above equations are accompanied by the initial and boundary conditions.
Particularly, we set $\w(x,0) = \w_0(x)$ in $\Om$,
where $\w_0$ is the given function, typically as a perturbation in terms of the potential
temperature of a steady-state flow where the mean values $\bar{\tempP}$ and $\bar{\pressE}$
of the Exner pressure and the potential temperature, respectively, 
are in a hydrostatic balance, cf.~\cite{GiraldoRostelli_JCP08}, i.e.,
\begin{align}
  \label{balance}
  \ccP \bar{\tempP} \frac{\dd \bar{\pressE}}{\dd x_2} = - g.
\end{align}

Moreover, we consider several types of boundary conditions on the boundary $\gom:=\pd\Om$:
\begin{itemize}
\item the {\em no-flux} or {\em reflecting boundary condition} $\bkv\cdot\bkn=0$,
  where $\bkn$ is the outer normal to the boundary,
\item the {\em nonreﬂecting boundary conditions} representing an 
  artificial boundary of an infinite domain,
\item the {\em periodic boundary conditions},
\item the {\em symmetric boundary conditions}.
\end{itemize}
Their practical realization in the context of the discontinuous Galerkin method
is described in Section~\ref{sec:semi}.

\section{Discontinuous Galerkin discretization}
\label{sec:DGM}

We solve numerically the system of equations \eqref{eq:NS} --  \eqref{eq:tau} by the
{\em space-time discontinuous Galerkin method}.
We refer to \cite[Chapters~8--9]{DGM-book}
for the detailed discretization of the Navier-Stokes equations
without the external gravity forces.
Here,  only the final necessary formulas are presented.

\subsection{Functional spaces}
Let $0=t_0<t_1<\ldots <t_r=T$ be a partition of $(0,T)$. We set 
$I_m =(t_{m-1},t_{m})$ and  $\tau_m = t_m- t_{m-1}$ for  $m=1,\ldots,r$.
For each $t_m,\ m=0,\ldots,r$, we consider a 
space partition $\Thm$ consisting of a finite number 
of closed 
triangles $K$ with mutually disjoint interiors and covering $\oO$,
the meshes $\Thm$ can vary for $m=0,\dots,r$, in general.
By $\bknK$ we denote the unit outer normal to the element boundary $\dK$, $K\in\Thm$.

For each $K\in\Thm$, $m=0,\dots,r$, we assign a positive integer number $p_K$
denoting the polynomial approximation degree on $K$. Then we define the spaces of
discontinuous piecewise-polynomial functions by
\begin{align}\label{Smhp}
  \Smhp = \{\vas: \Omega\to \R;\ \vas(x)|_K\in P_{\pK}(K)\ \forall\,K\in\Thm\},
  \quad \bSmhp=[\Smhp]^4,\qquad m=0,\dots,r,
\end{align}
where $P_{\pK}(K)$ denotes the space of all polynomials on $K$ of  degree $\leq \pK$.

The set of all edges of elements of mesh $\Thm$ is denoted by  $\Fhm$, 
the symbols $\FhmI$ and $\FhmB$ denote its subsets consisting
of all interior and boundary edges of $\Thm$, respectively.
Moreover, $\bpsi|_{\dK}\plus$ and $\bpsi|_{\dK}\minus$ are
the interior and exterior traces of $\bpsi\in\bSmhp$ on $\dK,\ K\in\Thm$, respectively. 
The symbols $\aver{\bpsi}$ and $\jump{\bpsi}$ denote
the mean value and the jump of functions $\bpsi\in\bSmhp$ on edges
$\gamma\subset \dK,\ K\in\Thm$ given by 
\begin{align}
  \label{jump}
  \begin{array}{llll}
    &\aver{\bpsi}_\gamma = \tfrac12\left(\bpsi|_{\dK}\plus + \bpsi|_{\dK}\minus\right),\ 
    &\jump{\bpsi}_\gamma = \left(\bpsi|_{\dK}\plus - \bpsi|_{\dK}\minus\right)\bknK\
    &\mbox{for }\gamma\in\FhmI, \\[2mm]
    &\aver{\bpsi}_\gamma = \bpsi|_{\dK}\plus,\
    &\jump{\bpsi}_\gamma = \bpsi|_{\dK}\plus \bknK\ 
    &\mbox{for } \gamma\in\FhmB.\\     
  \end{array}
\end{align}

Furthermore, let $q\ge0$ be the polynomial approximation degree with respect to time.
We define the space of discontinuous piecewise polynomial
functions on the space-time domain $Q_T$ by
\begin{align}
  \label{Shpq}
  \bShpq = & \left\{ \bpsi:\Om\times (0,T) \to \R^4;\
  \bpsi(x,t)|_{\Om\times I_m} = \sum\nolimits_{l=0}^q  t^l\, \bvas_k(x),\ \ \bvas_k\in \bSmhp,
  \ l=0,\dots, q,  \ m=0,\dots,r \right\}.
\end{align}
Since $\bpsi\in \bShpq$ is discontinuous with respect to the time variable at
$t_m,\ m=1,\dots, r-1$, we define the jump  with respect to time
on the time level $t_m,\ m=0,\dots, r-1$ by
\begin{align}
  \label{tjump}
  \tjump{\bpsi}_m := \bpsi|_m^+ - \bpsi|_m^-,
  \qquad \bpsi|_m^\pm := \lim_{\delta\rightarrow 0\pm}\bpsi(t_m+\delta),
\end{align}
where $\bpsi|_0^-$ is typically taken from the initial condition.

\subsection{Space semi-discretization}
\label{sec:semi}

The space semi-discretization of the system \eqref{eq:NS} employs the following
properties and relations: 
\begin{itemize}
\item the convective fluxes fulfill $\bF_i (\w) = \bA_i(\w) \w$,  where $\bA_i$ are the
  Jacobian matrices of $\bF_i (\w)$, $i=1,2$,
\item the physical inviscid flux through element boundary $\dK$, $K\in\Thm$
  having the normal vector $\bknK=(n_1,n_2)$
  is approximated by the  numerical flux
  \begin{align}
    \label{num_flux}
    \sum\nolimits_{i=1}^2 \bF_i(\w) n_i|_{\dK}
  \approx\PP\Plus(\w,\bknK)\w|_\gamma\plus + \PP\Minus(\w,\bknK)\w|_\gamma\minus,
  \end{align}
  where
  $\PP\PMinus$ are matrices from a suitable
  decomposition of the matrix $\sum_{i=1}^2 \bA_i(\w) n_i$, e.g.,
  Lax-Friedrichs or Vijayasundaram numerical fluxes can be written in this form \cite{FEI2},
\item
  the diffusive terms in \eqref{eq:NS} are discretized by the symmetric interior penalty
  variant of DG  method.
\end{itemize}

By symbol $\LSP{\bku}{\bpsi}{M}$ we denote the $L^2(M)$-scalar product of vector-valued functions
$\bku,\bpsi: Q_T\to \R^4$ over one- or two-dimensional set $M$.
In order to simplify the notation, we introduce
the symbol $\bAv=(\bA_1, \bA_2)$ denoting the ``vector of matrices'' (third-order tensor).
Then
\begin{align}
  \LSP{\bAv(\w) \w }{\nabla \bpsi}{M} = \sum\nolimits_{i=1}^2\LSP{\bA_i(\w)\w}{ \pdd{\bpsi}{x_i}}{M}.
\end{align}
Similarly, in virtue of \eqref{Kij} and \eqref{jump}, we use the notation
\begin{align*}
  \LSP{\mKv(\w)\nabla \w }{\nabla \bpsi}{M}
  = \sum\nolimits_{i,j=1}^2\LSP{\mK_{i j}(\w) \pdd{\w}{x_j}}{ \pdd{\bpsi}{x_i}}{M}
  \ \ \mbox{and}
  \ \
  \LSP{\aver{\mKv(\bbw)\nabla \w}}{\jump{\bpsi}}{\gamma}
  = \sum\nolimits_{i,j=1}^2\LSP{\mK_{i j}(\bbw) \pdd{\w}{x_j}}{
    (\bpsi|_{\gamma}\plus - \bpsi|_{\gamma}\minus)}{\gamma}n_i,
\end{align*}
where $(n_1,n_2)$ is the outer normal to an edge $\gamma\subset \dK$.

We define the forms representing the space semi-discretization of the convective,
diffusive and gravity terms in \eqref{eq:NS} by discontinuous Galerkin method as
\begin{align}
  \label{ah}
  & \ahmL(\bbw, \w, \bpsi) :=  - \sumK{\bAv (\bbw) \w}{\nabla \bpsi} 
   +\sumk\left(\left({\PP\Plus(\bbw,\bknK)\w|_{\dK}\plus},{\bpsi}\right)_{\pd K}
  + \left({\PP\Minus(\bbw,\bknK)\w|_{\dK}\minus },{\bpsi}\right)_{\pd K\setminus\gom}\right)
  - \sumK{\bB\w}{\bpsi}
  \notag \\
  & \qquad + \sumK{\mKv(\bbw)\nabla \w}{\nabla \bpsi} 
  -\sumF{\aver{\mKv(\bbw)\nabla \w}}{\jump{\bpsi}}  
  -\sumF{\aver{\mKv(\bbw)\nabla \bpsi}}{\jump{\w}} 
   + \sumF{\sigma\jump{\w}}{\jump{\bpsi}}, 
    \\
   & \ahmt(\bbw,\bpsi) := 
   - \sumKK{ \PP\Minus(\bbw,\bknK) \w_D}{\bpsi}{\pd K\cap\gom}
   - \sumFB{ {\mKv(\bbw)\nabla \bpsi}}{{\w_D \bknK}}
   + \sumFB{\sigma \w_D}{\bpsi}, \notag \\
   & \ahm(\w, \bpsi) := \ahmL(\w, \w, \bpsi) -  \ahmt(\w,\bpsi),\qquad \bbw, \w,\bpsi\in\bSmhp,
   \quad m=1,\dots,r,
   \notag
\end{align}
where
$\sigma\sim p_K^2/ \mbox{diam}(K)$ is the penalty parameter
and $\w_D$ is the state vector following from the boundary conditions. 
We note that the form $\ahmL$ from \eqref{ah}
is linear with respect to its second and third arguments,
which is employed in the iterative solution of the arising nonlinear algebraic systems in
Section~\ref{sec:sol}.

The boundary state vector $\w_D$ is set in the following way.
\begin{itemize}
\item {\em No-flux boundary} conditions: the vector $\w_D$
  has the same density, energy and the tangential component of the velocity
  as $\w_{\dK}$ (= the trace of $\w$ on $\dK\subset \gom$)
  but the opposite normal component of the velocity, i.e., we set
  \begin{align}
    \label{BCflux}
    \w_D|_{\dK} =\Mir(\w|_{\dK},\bknK),
    \quad \mbox{where}\quad \Mir(\w,\bkn) =
    \left(\rho, \rho (\vv -2(\vv\cdot\bkn)\bkn), \ener\right)^{\TT}
    \quad \mbox{for}\quad
    \w=(\rho, \rho\vv, \ener)^{\TT}.
  \end{align}
\item {\em Periodic boundary conditions}:
  we set $\w_D|_{\gamma^+} = \w|_{\gamma^-}$ and $\w_D|_{\gamma^-} = \w|_{\gamma^+}$,
  where $\gamma^+$ and $\gamma^-$ denote the
  pair of the corresponding periodic edges.
\item {\em Symmetric boundary conditions}: $\w_D$ is set similarly as for the
  periodic boundary conditions, depending on the line of symmetry. E.g., if the symmetry
  is prescribed along the axes $x_1=0$, the vector $\w_D$
  has the same density, energy and the vertical component of the velocity $v_2$
  but the opposite horizontal component of the velocity $v_1$.
\item {\em Non-reflecting boundary conditions}: we employ the approach
  from \cite[Section~8.3.2]{DGM-book}, where $\w_D$ is the solution of the
  nonlinear Riemann problem $\BC^{\rm RP}(\w_\gom, \w_{\mathrm{init}})$
  defined for the Euler equations ($\bRv=0$ in \eqref{eq:NS}) without
  the source (gravity) terms. Here, $\w_{\dK}$ is the trace of $\w$
  on $\dK\subset \gom$ and $\w_{\mathrm{init}}$ is typically given by
  the initial (unperturbed) flow field.
\end{itemize}

We are aware that the aforementioned realization of the
non-reflecting boundary conditions is rather heuristic since the
treated Riemann problem does not involve the gravity term.
However, this approach works very well provided that
high polynomial approximation degrees are employed, which is documented by
the numerical experiments in Section~\ref{sec:BC}.

\subsection{Performance of the non-reflecting boundary conditions}
\label{sec:BC}

We consider
the steady-state solution of hydrostatic flow in 
the computational
domain $\Om=(-25\,600, 25\,600)\times (0, 6\,400)\si{m}$ with the constant potential temperature
$\bar{\tempP}= 300\si{K}$ and the velocity $\bkv = 0$.
The non-reflecting boundary conditions are prescribed on the vertical parts of the boundary domain
and the no-flux boundary conditions are treated on the horizontal ones.
We carried out computation
till $T=10\si{s}$ and investigated the ability of the proposed numerical scheme to keep
the initial steady-state solution. 
We employed $P_1$, $P_3$ and $P_5$ polynomial approximation with respect to space
using grids having $\#\Th=117$, 406, 2\,728 and 10\,916 triangles.
In Table~\ref{tab:BC}, we show the quantities
\begin{align}
  \label{num:BC}
  \min{\delta\tempP} = \min_{x\in\Om} \left(\tempP(x) - \bar{\tempP}\right)/\bar{\tempP}, \qquad
  \max{\delta\tempP} = \max_{x\in\Om} \left(\tempP(x) - \bar{\tempP}\right)/\bar{\tempP}, \qquad
  \Delta{\delta\tempP} = \max{\delta\tempP}  - \min{\delta\tempP}.
\end{align}

We observe the large violation of the constant flow field for $P_1$ approximation even for
relatively very fine grids. The decrease of $\Delta{\delta\tempP}$ with respect to the mesh
size step $h$ is slow. On the other hand, 
$P_5$ approximation gives almost negligible values of
$\Delta{\delta\tempP}$ even for the coarsest grid. This effect underlines the importance of the
use of higher order approximation and it nicely correlates with the used mesh adaptive algorithm,
described in Section~\ref{sec:AMA} which generates typically meshes having
large elements with high polynomial
approximation degrees for unperturbed (smooth) flows.
\begin{table}
  \begin{center}
    \setlength{\tabcolsep}{12pt}
    {\footnotesize
    \begin{tabular}{crrccc}
    \hline
    degree & $\#\Th$\  & $\DoF$\  & $\min{\delta \tempP}$ & $\max{\delta \tempP}$ & $\Delta({\delta \tempP})$ \\ 
    \hline
    $P_1$ & 117 & \np{351} & -1.0672E-02 & 3.9217E-05 & 1.0711E-02 \\
    $P_1$ & 703 & \np{2109} & -2.2365E-03 & 3.3900E-05 & 2.2704E-03 \\
    $P_1$ & \np{2728} & \np{8184} & -4.6223E-04 & 1.3397E-06 & 4.6357E-04 \\
    $P_1$ & \np{10916} & \np{32748} & -1.2871E-04 & 1.1181E-06 & 1.2983E-04 \\
    \hline           
    $P_3$ & 117 & \np{1170} & -3.3737E-06 & 4.8027E-07 & 3.8539E-06 \\
    $P_3$ & 703 & \np{7030} & -1.6309E-07 & 2.6487E-08 & 1.8957E-07 \\
    $P_3$ & \np{2728} & \np{27280} & -8.3347E-09 & 1.2776E-09 & 9.6122E-09 \\
    \hline           
    $P_5$ & 117 & \np{2457} & -2.1223E-10 & 1.7873E-11 & 2.3010E-10 \\
    $P_5$ & 703 & \np{14763} & -2.6707E-12 & 6.0177E-13 & 3.2725E-12 \\
    \hline
    \end{tabular}
    }
  \end{center}
  \caption{Performance of the non-reflecting boundary conditions, the violence of the
    constant potential temperature field, cf.~\eqref{num:BC}. The value $\DoF$ means the number
    of degrees of freedom of the space semi-discretization per one equation.}
\label{tab:BC}
\end{table}

\subsection{Full space-time discretization}

To define the full space-time DG discretization of \eqref{eq:NS}, we introduce the
forms
\begin{align}
  \label{Ah}
   \Ahm(\w, \bpsi) & :=  \int_{I_m}\Bigl( \Lsp{\pdt{\w}}{\bpsi} + \ahm(\w,\psi)\Bigr)\dt
  + \LSP{\tjump{\w}_{m-1}}{\bpsi|_{m-1}^{+}}{\Om}, \\
  \AhmL(\bbw, \w, \bpsi) & :=  \int_{I_m}\left( \Lsp{\pdt{\w}}{\bpsi}
  + \ahmL(\bbw,\w,\psi)\right)\dt
   + \LSP{\tjump{\w}_{m-1}}{\bpsi|_{m-1}^{+}}{\Om}, \quad \bbw,\w,\bpsi\in\bShpq,\quad
   m=1,\dots,r.\notag
\end{align}
We note that form $\AhmL$ 
is linear with respect to its second and third arguments.

\begin{definition}
  We say that $\wht \in \bShpq$ is the {\em space-time discontinuous Galerkin} solution
  of \eqref{eq:NS} if
  \begin{align}
    \label{STDGM}
     \Ahm(\wht,\bpsi_h) = 0 \qquad \forall \bpsi_h \in\bShpq, \quad m=1,\dots,r,
  \end{align}
  where we set $\wht|_0^- := \w_0$ (= the initial condition).
\end{definition}

\begin{figure}
  \includegraphics[width=0.490\textwidth]{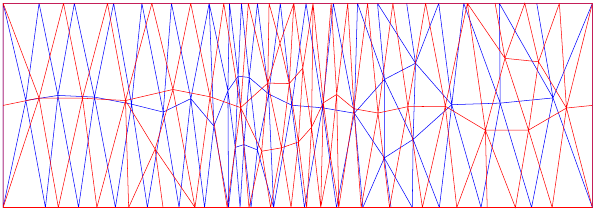}
  \includegraphics[width=0.490\textwidth]{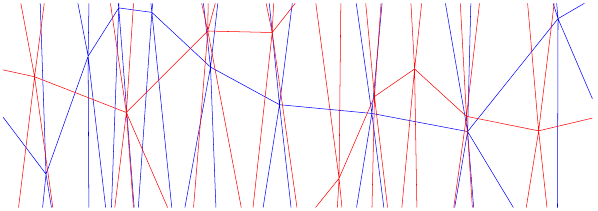}
  \caption{Illustration of meshes $\ThmM$ (red) and $\Thm$ (blue),
    total view (left) and a detail (right).}
  \label{fig:IGW_meshes}
\end{figure}

The term $\LSP{\tjump{\wht}_{m-1}}{\bpsi|_{m-1}^{+}}{\Om}$ in \eqref{Ah}
represents the time penalty which
joins together the approximate solution on two subsequent time layers $I_{m-1}$ and $I_m$
in a weak sense.
This term is applicable also for varying spaces $\bSmhpM$ and $\bSmhp$.
Indeed, due to \eqref{tjump}, we have
\begin{align}
  \label{Tjump}
  \LSP{\tjump{\wht}_{m-1}}{\bpsi|_{m-1}^{+}}{\Om} =
  \LSP{ \wht|_{m-1}^+}{\bpsi|_{m-1}^{+}}{\Om} -\LSP{ \wht|_{m-1}^-}{\bpsi|_{m-1}^{+}}{\Om}. 
\end{align}
If  $\bSmhpM \not= \bSmhp$ then the last term in \eqref{Tjump} exhibits the
integration of the product of two functions which are piecewise polynomial on different meshes
$\ThmM$ and $\Thm$.
In Section~\ref{sec:AMA}, we introduce the {\hpAMA} algorithm which generates completely
different meshes for two sequential time steps, see 
Figure~\ref{fig:IGW_meshes} for an illustration.
In this case, the evaluation of the last integral in \eqref{Tjump} is a little delicate
since we integrate $\wht|_{m-1}^-$ over $K\in\Thm$ but  $\wht|_{m-1}^-$ is discontinuous on $K$.
In \cite{AMAtdp}, we presented a strategy, which employs a non-adaptive
composite quadrature rule. The evaluation of this integral exhibits another source
of computational errors. However, the numerical experiments in
Section~\ref{sec:quantity} show that the inaccuracy of the mass and energy conservation
is on the same level as the computation on a fixed mesh. On the other hand, if
$\ThmM = \Thm$, the numerical integration 
of the terms in \eqref{Tjump} can be carried out exactly (up to rounding errors)
since the integrated functions are polynomial.

The STDG discretization \eqref{STDGM} represents the system of $\Nhm$
algebraic equations for each $m=1,\dots,r$, where
\begin{align}
  \label{Nhm}
  \Nhm = (q+1)\Mhm\quad\mbox{with}\quad  \Mhm =  \sum_{K\in\Thm}\nolimits (\pK+1)(\pK+2)/2.
\end{align}
On the contrary,
the diagonally implicit Runge-Kutta methods give  $q$ systems of size $\Mhm$ 
and multi-step methods give only one system of size $\Mhm$, hence their solution is cheaper.
This is the main disadvantage of the time discontinuous Galerkin discretization.
We note that if $q=0$ (piecewise constant in time) then the time DG scheme is identical
to the backward (implicit) Euler method. 

On the other hand, the STDG method has the super-convergence property in the time
nodes $t_m,\ m=1,\dots,r$, namely 
$O(\tau^{2q+1})$. Therefore, the use
of piecewise linear approximation ($q=1$)
seems to be sufficient for the majority of numerical experiments since
it provides a third-order method at time nodes. 
Moreover, whereas multi-step of Runge-Kutta methods give the approximate solution only at
time nodes $t_m,\ m=1,\dots,r$, the STDG discretization
provides an approximate solution on the whole interval $I_m$, $m=1,\dots,r$ directly.
It is advantageous in the definition of the space and time discretization errors,
cf.~Section~\ref{sec:etas}.

\section{Solution strategy}
\label{sec:sol}
In this section, we describe the strategy for solving the discrete problem \eqref{STDGM}.
We omit the details and we refer to \cite{stdgm_est}, \cite[Chapters~8-9]{DGM-book} for
a complete description.

\subsection{Solution of nonlinear algebraic systems}

Let $m=1,\dots, r$ be an arbitrary but fixed index of the time level.
The nonlinear algebraic system \eqref{STDGM} 
is solved iteratively by a 
Newton-like method where we do not differentiate the form $\Ahm$ but employ
its linearization $\AhmL$, cf. \eqref{Ah}. Particularly, let $\wm:=\wht|_{\Om\times I_m}$
denote the restriction of $\wht$ on the space-time layer $\Om\times I_m$.
Then we define the sequence $\wmk$, $k=0,1,\dots,$ approximating $\wm$ by
\begin{subequations}
  \label{SS1}
  \begin{align}
  \label{SS1a}
  \wmN &= \Ext(\wmNN),  \\
  \label{SS1b}
  \wmk &= \wmkM + \lambda_k \buk ,\qquad k=1,2,\dots \\
  \label{SS1c}
  &\quad \mbox{where } \buk\in\bShpq:\ \ 
  \AhmL(\wmkM, \buk, \bpsi_h) = - \Ahm(\wmkM, \bpsi_h)\quad  \forall \bpsi_h \in\bShpq. 
  \end{align}
\end{subequations}
The symbol $\Ext$ in formula \eqref{SS1a} means that
the initial guess $\wmN$ is extrapolated either
from the previous time layer or from the initial condition.
If $\bSmhpM \not= \bSmhp$, we have to project $\wmNN\in \bSmhpM$
to the new space $\bSmhp$ but we do not emphasize it explicitly.
Moreover, the function $\buk$ in \eqref{SS1b} is the update in the $k$-th iteration
and it is obtained from the linear algebraic system \eqref{SS1c}.
The convergence is improved by the
suitable choice of the damping factor $\lambda_k \in (0,1]$. Its choice is
  controlled by the monitoring factor \cite{newton}
  \begin{align}
    \zeta_k = \frac{ \norm{\left\{\Ahm(\wmk,\bpsi_i)\right\}_{i=1}^{\Nhm}}{\ell^2} }
    { \norm{\left\{\Ahm(\wmkM,\bpsi_i)\right\}_{i=1}^{\Nhm}}{\ell^2} },
  \end{align}
  where $\bpsi_i$, $i=1,\dots,\Nhm$ are the basis functions on $\Om\times I_m$ such that
  $\norm{\bpsi_i}{L^2(\Om\times I_m)} = 1$ and $\norm{\cdot}{\ell^2}$ is the Euclidean
  discrete norm.  Hence, the damping factor $\lambda_k$ has to be chosen such that
  $\zeta_k < 1$ for each iterates $k = 1,2,\dots$.

  It is advantageous to solve the linear algebraic system  \eqref{SS1c} iteratively
since $\buk=0$  is a suitable initial guess and only a few iterative steps are necessary.
We employ the  GMRES method with block ILU(0) preconditioner but other solvers are possible.
The stopping criterion for the iterative process
\eqref{SS1} is presented in Section~\ref{sec:etas}.
We note that it is not necessary to update the matrix on the left-hand side of
\eqref{SS1c} for each $k$. It means that in the first argument of $\AhmL$, we can employ
a former approximation $\wmkMs$ for some $s>1$. Particularly, if the monitoring
function $\zeta_k$ is sufficiently smaller than 1, we do not update the matrix on the
left-hand side of \eqref{SS1c}.

The presented approach can be treated also in the context of the IMEX formulation.
Let $\wht$ be the approximate solution given by \eqref{STDGM} and $\lambda\in(0,1]$ then
we have the identity
\begin{align}
  \label{IMAX1}
   \AhmL(\wht,\wht,\bpsi_h) =  \AhmL(\wht,\wht,\bpsi_h) -\lambda \Ahm(\wht,\bpsi_h) 
  \qquad \forall \bpsi_h \in\bShpq, \quad m=1,\dots,r.
\end{align}
Treating $\wht$ in the second argument of $\AhmL$ on the left-hand side of \eqref{IMAX1} implicitly
and in the other arguments explicitly by extrapolating
from the previous time layer (cf.~\eqref{SS1a}),
we obtain the linear systems
\begin{align}
  \label{IMAX2}
  \AhmL(\Ext(\wmNN),\wm,\bpsi_h) =  \AhmL(\Ext(\wmNN),\Ext(\wmNN),\bpsi_h)
  -\lambda \Ahm(\Ext(\wmNN),\bpsi_h) 
  \qquad \forall \bpsi_h \in\bShpq, \quad m=1,\dots,r.
\end{align}
The formulation \eqref{IMAX2} is equivalent to the first iterative step of method \eqref{SS1} for
$k=1$.
So the additional iterates in  \eqref{SS1} for $k\ge 2$ exhibits an improvement of the accuracy
which can be controlled by the stopping criterion presented in Section~\ref{sec:etas}.
In our experience, the fully-implicit technique \eqref{SS1} is more stable than the
semi-implicit one \eqref{IMAX2}, particularly it is possible to use larger time steps.
We demonstrate this effect by a numerical example in Section~\ref{sec:stabil}.

\subsection{Choice of the time step and iterative stopping criterion}
\label{sec:etas}
The accuracy and efficiency of the whole computational process strongly depend on
the choice of the time step and the stopping criterion for the iterative solver
\eqref{SS1}. A too large time step and/or a too weak stopping criterion cause a loss of accuracy.
On the other hand, a too small time step and/or a too strong  stopping criterion lead to an
increase in the computational time without a significant impact on the accuracy.
In order to balance these two aspects, we employ residual
based estimators indicating the errors arising from the spatial and temporal discretizations
and the inaccuracy of iterative algebraic solvers.

For each time interval $I_m,\ m=1,\dots,r$, we define the estimators
\begin{align}
  \label{etas}
  \etaAm(\wht) = \sup_{0\not=\bpsi_h\in \bShpq}
  \frac{\Ahm(\wht,\bpsi_h)}{\norm{\bpsi_h}{X}}, \ \ 
  \etaSm(\wht) = \sup_{0\not=\bpsi_h\in \bShpqS}
  \frac{\Ahm(\wht,\bpsi_h)}{\norm{\bpsi_h}{X}}, \ \
  \etaTm(\wht) = \sup_{0\not=\bpsi_h\in \bShpqT}
  \frac{\Ahm(\wht,\bpsi_h)}{\norm{\bpsi_h}{X}}, 
\end{align}
where  the space of functions $\bShpq$ is given by \eqref{Shpq} and
the spaces $\bShpqS$ and $\bShpqT$  are given analogously by an increase of the polynomial
approximation degree by one with respect to space and time, respectively. Particularly, we have
\begin{align}
  \label{Shppqq}
  \bShpqS = & \left\{ \bpsi:\Om\times (0,T) \to \R^4;\
  \bpsi(x,t)|_{\Om\times I_m} = \sum\nolimits_{l=0}^q  t^l\, \bvas_k(x),\ \ \bvas_k\in \bSmhpP,
  \ l=0,\dots, q,  \ m=0,\dots,r \right\}, \\
  \bShpqT = & \left\{ \bpsi:\Om\times (0,T) \to \R^4;\
  \bpsi(x,t)|_{\Om\times I_m} = \sum\nolimits_{l=0}^{q+1}  t^l\, \bvas_k(x),\ \ \bvas_k\in \bSmhp,
  \ l=0,\dots, q+1,  \ m=0,\dots,r \right\}, \notag
\end{align}
where (cf.~\eqref{Smhp})
\begin{align}\label{SmhpP}
  \bSmhpP=[\SmhpP]^4\quad\mbox{with}\quad
  \SmhpP = \{\vas: \Omega\to \R;\ \vas(x)|_K\in P_{\pK+1}(K)\ \forall\,K\in\Thm\},
  \qquad m=0,\dots,r.
\end{align}

 Moreover, the norm $\norm{\cdot}{X}$ is defined by 
\begin{align}
  \label{normX}
  \norm{\bpsi_h}{X}=\left(
  \normP{\bpsi_h}{L^2(\Om\times I_m)}{2} +
  \normP{\nabla \bpsi_h}{L^2(\Om\times I_m)}{2} +
  \normP{\pdt\bpsi_h}{L^2(\Om\times I_m)}{2}\right)^{1/2}.
\end{align}
The quantity  $\etaAm(\wht)$ represents the
inaccuracy of the solution of the algebraic system \eqref{STDGM} since it vanishes
if  \eqref{STDGM}  is solved exactly. In the same spirit, $\etaSm(\wht)$ and $\etaTm(\wht)$
represent the inaccuracy of the spatial and temporal discretization, respectively.
Particularly, they vanish if the approximate solution  is equal to the exact one.
The numerical study of these error estimators was presented in \cite{stdgm_est}.
We note that due to the choice of the norm \eqref{normX}, the estimators in \eqref{etas}
can be evaluated as the sum  of local estimators computed for each space-time
element $K\times I_m$, $K\in\Thm$ separately.
These local estimators are computable by standard tools of the
constrained minimization since the corresponding systems of equations are small and
thus the evaluation of \eqref{etas} is fast.

Employing   estimators \eqref{etas}, we define the stopping criterion of the iterative method
\eqref{SS1} by
\begin{align}
  \label{CA}
  \etaAm(\wht) \leq c_A \min\left(\etaSm(\wht), \etaTm(\wht)\right),\quad m=1,\dots, r,
\end{align}
where $c_A\in(0,1)$.
It means that the algebraic error estimator has to be bounded with respect 
to the spatial and temporal ones.
Similarly, we choose adaptively the time step $\tau_m$, $m=1,\dots, r$ such that
spatial and temporal error estimators are balanced, particularly, we set the time step
such that
\begin{align}
  \label{CT}
  \etaTm(\wht) \approx c_T \etaSm(\wht),\quad m=1,\dots, r,
\end{align}
for some $c_T\in(0,1)$. In practice, we employ the values $c_A=0.01$ and $c_T=0.2$.

\subsection{Numerical study of the stability and adaptivity of the STDG scheme}
\label{sec:stabil}

\begin{figure}
  \includegraphics[width=0.470\textwidth]{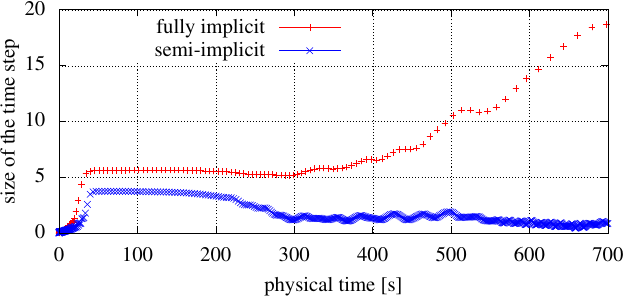}
  \hspace{0.02\textwidth}
  \includegraphics[width=0.470\textwidth]{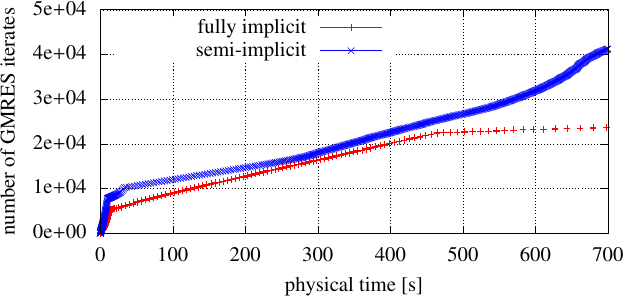}
  \caption{Comparison of the semi-implicit (blue) and fully implicit (red) schemes:
   the dependence of the size of the time step $\tau_m$ (left) and
   the accumulated number of GMRES iterations (right) with respect to the physical time
   $t\in(0,700)\si{s}$. Each dot corresponds to one time step $m=1,\dots, r$.}
  \label{fig:stab}
\end{figure}

To demonstrate the stability and the adaptive performance of the
fully implicit scheme, we consider the rising thermal bubble
from Section~\ref{sec:thermal} with the smooth initial condition \eqref{thermal1}.
We consider a fixed quasi-uniform grid having 556 elements, $P_2$ approximation
with respect to space and $P_1$ approximation with respect to time variables.
Therefore, the size of the algebraic system is $\Nhm= \np{26688}$.
We compare the computations obtained by
\begin{itemize}
\item {\em fully implicit scheme} \eqref{SS1} using the adaptive choice of the time step
  according to \eqref{CT} with $c_T= 0.2$ and adaptive stopping criterion \eqref{CA} with
  $c_A=0.01$,
\item {\em semi-implicit scheme} \eqref{IMAX2} using the adaptive choice of the time step
  according to \eqref{CT} with $c_T= 0.1$ which guarantees the stability and a comparable
  quality of the results as the fully implicit scheme -- the larger value of $c_T$ gives
  an unstable solution.
\end{itemize}
Figure~\ref{fig:stab} shows the size of the time step $\tau_m$ and
the accumulated number of GMRES iterations for each time step $m=1,2,\dots$.
We observe that the fully implicit scheme admits the use of a larger size of the time steps
and the corresponding number of GMRES iterations is smaller. We note that the size of
$\tau_m$ in the last time steps corresponds to $\mbox{CFL}= 400$ in the ``common''
stability condition \cite[Section~8.4]{DGM-book}
\begin{align}
  \tau_m \leq  \min_{\gamma \subset \dK} \frac{\mbox{CFL}\, |K|}{|\gamma|\, \Lambda(\wht,\bknK)}\qquad
  \forall \gamma\subset\dK\quad \forall K\in\Thm,
\end{align}
where
$|K|$ denotes the area of $K\in\Thm$, $|\gamma|$ is the length of an edge $\gamma\subset\dK$ and
$\Lambda(\wht,\bknK)$ is the maximal
eigenvalue of matrix $\PP=\PP\plus+\PP\minus$ arising in the numerical flux \eqref{num_flux}.
The computational time
for the semi-implicit method was $\np{2437}\si{s}$ whereas for the implicit one
$\np{1363}\si{s}$.

\begin{figure}
  \includegraphics[width=0.470\textwidth]{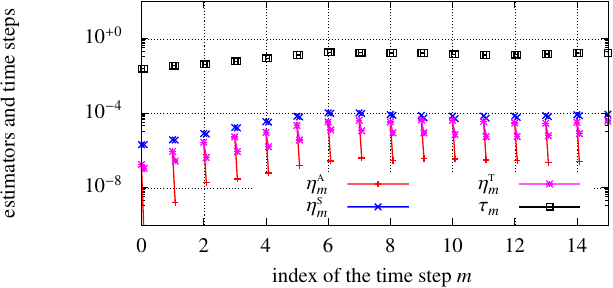}
  \hspace{0.02\textwidth}
  \includegraphics[width=0.470\textwidth]{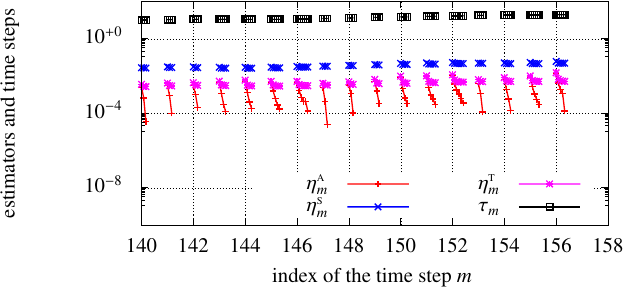}

  \caption{Performance of the adaptive algebraic stopping criteria
    and the adaptive choice of the time step: estimators $\etaAm$, $\etaSm$, $\etaTm$
    (cf.~\eqref{etas}) and  $\tau_m$ for $m=1,\dots, 15$ (left) and
    $m=141,\dots, 156$ (right), each dot corresponds to one iterative step $k$ in \eqref{SS1}.}
  \label{fig:stab1}
\end{figure}

Moreover, Figure~\ref{fig:stab1} shows the performance of the adaptive algebraic stopping criterion
\eqref{CA}
and the choice of the time step \eqref{CT}.
We plot the size of estimators $\etaAm$, $\etaSm$, $\etaSm$
(cf.~\eqref{etas}) and the size of the time step $\tau_m$ for several starting and final
time levels $m$. We observe that the size of the time steps and algebraic errors are controlled
by the aforementioned conditions.
Moreover, only 2-3 nonlinear iterations are performed
for smaller $m$ when the time steps are shorter and about $3-5$ steps at the end of the
computation where the time steps are larger.

\section{Mesh adaptation}
\label{sec:AMA}

Mesh adaptive methods try to reduce the number of
degrees of freedom ($\DoF$) while the accuracy is preserved. 
This is beneficial namely for the simulation of evolution problems with traveling phenomena
having sharp interfaces. The idea of the presented 
adaptation procedure is to generate the sequence
of the spaces $\Smhp$, $m=0,\dots,r$ (i.e., the meshes $\Thm$ and the corresponding
polynomial degrees $\pK$, $K\in\Thm$) such that
\begin{enumerate}[label=({\roman*})]
\item the interpolation error of a selected scalar quantity
  is under  the given tolerance $\TOL>0$, \label{c1}
\item the number of degrees of freedom ($=\Nhm$, cf. \eqref{Nhm}) is minimal for each $m=0,\dots,r$.
\end{enumerate}
More precisely, let $\qht=\qht(\wht):\Om\times(0,T)\to\R$ denote a scalar
solution-dependent quantity, e.g., the density and the Mach number are commonly
used for aerodynamic problems,
cf.~\cite{BelmeDervieuxAlauzet_JCP12,AlauzetLoseille_JCP18,FreyAlauzet,HabashiAll}.
Then condition \ref{c1} means
\begin{align}
  \label{AMA1}
  \norm{\Wqht(t) - \Pi \Wqht(t)}{L^2(\Om)} \leq \TOL \quad \forall t\in I_m,\ m=1,\dots,r,
\end{align}
where $\Wqht(t)\in\SmhpP$, cf.~\eqref{SmhpP}
is the higher-order reconstruction of the scalar piecewise polynomial function
$\qht(t)=\qht(\wht(\cdot,t))$ for $t\in I_m$ and $\Pi$ is the $L^2(\Om)$-projection to $\Smhp$. 
In practice, we relax condition \eqref{AMA1} and require
\begin{align}
  \label{AMA2}
  \norm{\Wqht(t) - \Pi \Wqht(t)}{L^2(\Om)} \leq \TOL \quad \forall t\in \cJ_m,\ m=1,\dots,r,
\end{align}
where $\cJ_m$ is the finite set of selected $t\in I_m$.
Typically $\cJ_m$ contains endpoints of $I_m$
and the time integration nodes. The reconstruction $\Wqht$
is carried out locally by a least-square technique, cf. \cite[Section~9.1.1]{AMA-book}.

The whole adaptive procedure is described by Algorithm~\ref{alg:hp-AMA-ST}.
We start with an initial mesh, perform one time step by solving problem \eqref{STDGM},
evaluate $\qht = \qht(\wht)$ 
and reconstruct $\Wqht(t)$, $t\in\cJ_m$.
If condition \eqref{AMA2} is valid we proceed to the next time step
using the same mesh and polynomial approximation degrees. Otherwise, we carry out
the re-meshing using  the {\em anisotropic $hp$-mesh adaptation} ($hp$-AMA)
technique (step \ref{alg:G1x} of Algorithm~\ref{alg:hp-AMA-ST}),
whose theoretical as well as practical issues
are described in detail in the recent
monographs \cite{AMA-book}, see also \cite{AMAtdp}.

\begin{algorithm}[ht]
  \caption{Space-time anisotropic $hp$-mesh adaptive algorithm.}
  \label{alg:hp-AMA-ST}
\begin{spacing}{1.15}
  \begin{algorithmic}[1]
    \STATE inputs: tolerance $\TOL>0$, initial time step $\tau_1$, initial
    $hp$-mesh $\Thn$ and $\pK\ge1$  $\forall K\in\Thn$
    \STATE set $m=1$, $t=0$ 
    \WHILE {$t < T$}
    \REPEAT \label{alg:G1a}
    \REPEAT \label{alg:G1b}
    \STATE 
    solve $\Ahm(\wht, \bpsi_h)= 0$\ \ $\forall\bpsi_h\in\bShpq$,
    cf. \eqref{STDGM} until \eqref{CA} is fulfilled \label{alg:G1c}
    \STATE evaluate $\etaSm(\wht)$ and $\etaTm(\wht)$ by \eqref{etas}
    \IF {$\etaTm > c_T \etaSm$}
    \STATE reduce the size of the time step $\tau_m$
    \ENDIF
    \UNTIL{$\etaTm \leq c_T \etaSm$}
  \STATE evaluate $\Wqht$ by a higher-order reconstruction of the quantity $\qht(\wht)$
  \STATE set $\etam:= \max_{t\in\cJ_m}\norm{\Wqht(t) -\Pi\Wqht(t)}{L^2(\Om)}$\label{alg:G1w}
  \IF{$\etam \leq \TOL$} \label{alg:G1z}
  \STATE set $t:= t+\tau_m$, propose the new size of the time step $\tau_{m+1}$
  \STATE
  set $\ThmP:= \Thm$ and the same $p_K,\ K\in\Thm$, set $m:= m+1$
  \ELSE
  \STATE generate a new mesh $\Thm$ \& new polynomial degrees  $p_K,\ K\in\Thm$
  by $hp$-AMA technique\label{alg:G1x}
  \STATE reduce the size of the time step $\tau_m$
  \ENDIF
  \UNTIL{$\etam \leq \TOL$}
\ENDWHILE
\end{algorithmic}
\end{spacing}
\end{algorithm}

The $hp$-AMA method sets an optimal size of mesh elements by the equidistribution
of an interpolation error estimate, i.e., the aim is to fulfill the condition
\begin{align}
  \label{AMA1a}
  \max_{t\in\cJ_m}\norm{\Wqht(t) - \Pi \Wqht(t)}{L^2(K)} \approx \frac{\TOL}{\sqrt{\#\Th}} \quad
    \forall K\in\Thm,\ m=1,\dots,r,
\end{align}
where $\#\Th$ is the number of elements of $\Thm$.
Moreover, the shape and orientation of mesh elements are optimized locally in such a way that
for each $K\in\Thm$, we fix its size and seek a new shape and orientation by minimizing
the interpolation error estimate. Similarly, we can modify the polynomial approximation
degree such that we fix the density of the number of degrees of freedom and
minimize the interpolation error estimate. This estimate depends on the
($\pK+1$)-directional derivatives of the reconstructed solution $\Wqht$.

There is a natural question of which quantity $\qht$ should be used for the higher-order
reconstruction.  
Solving non-hydrostatic mesoscale atmospheric problems, the very frequent quantity of
interest is the perturbation of the potential temperature $\tempP$, cf.~\eqref{eq:tempP}
with respect to the initial state. Hence, it seems to be favorable to employ $\tempP$
in the {\hpAMA} algorithm. Nevertheless, numerical experiments  in
Section~\ref{sec:quantity} demonstrate that the simple choice $\qht$ as the density
is also possible.

\begin{remark}
  We are aware that the computational time necessary for the solution of algebraic systems
  (the most time-consuming part of the whole computational process) does not
  depend only on the size of these systems. It is influenced also by the
  shape of elements, mesh grading, polynomial approximation degrees and the size of time steps.
  Nevertheless, the reduction of the number of $\DoF$ typically accelerates numerical simulations.
\end{remark}

\subsection{Numerical study of the mass and energy conservation and
  the choice of the quantity of interest}
\label{sec:quantity}

\begin{figure}
  \begin{center}

    \vertical{\hspace{15mm} mesh}
  \includegraphics[width=0.23\textwidth]{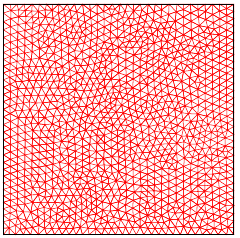}
  \includegraphics[width=0.23\textwidth]{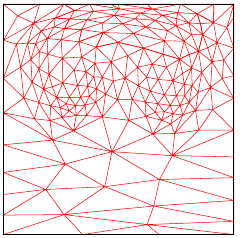}
  \includegraphics[width=0.23\textwidth]{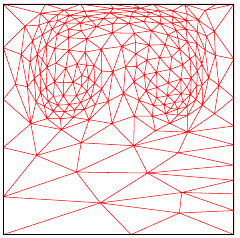}
  \includegraphics[width=0.23\textwidth]{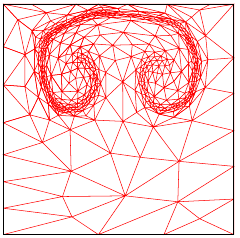}
  
  \vertical{\hspace{2mm} potential temperature}
  \includegraphics[width=0.23\textwidth]{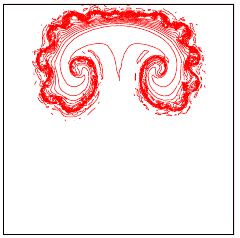}
  \includegraphics[width=0.23\textwidth]{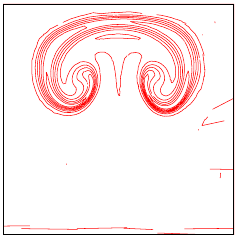}
  \includegraphics[width=0.23\textwidth]{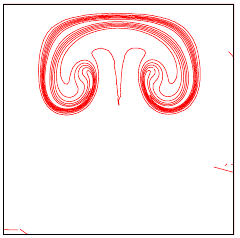}
  \includegraphics[width=0.23\textwidth]{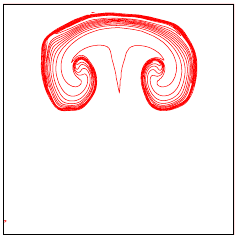}

  \vertical{\hspace{15mm} density}
  \includegraphics[width=0.23\textwidth]{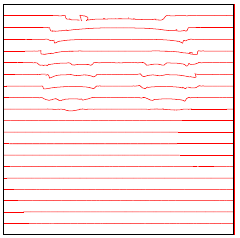}
  \includegraphics[width=0.23\textwidth]{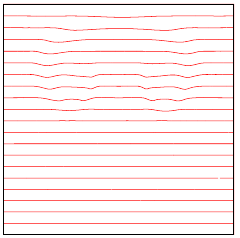}
  \includegraphics[width=0.23\textwidth]{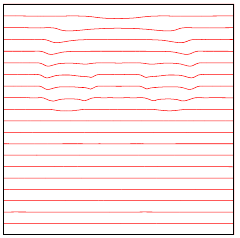}
  \includegraphics[width=0.23\textwidth]{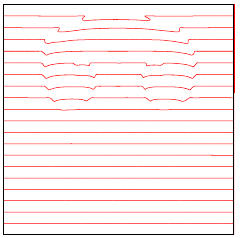}

  {\footnotesize
  \setlength{\tabcolsep}{14pt}
  \hspace{2mm}
  \begin{tabular}{llll}
  uniform  &
  \hpAMA, $\TOL=0.02$  & 
  \hpAMA, $\TOL=0.015$  & 
  \hpAMA, $\TOL=0.01$
  \\[0mm]
  $\Delta_{\mathrm{mass}}$ =  9.761E-06 & 
  $\Delta_{\mathrm{mass}}$ =  9.425E-06 & 
  $\Delta_{\mathrm{mass}}$ =  9.537E-06 & 
  $\Delta_{\mathrm{mass}}$ =  9.649E-06  
  \\
  $\Delta_{\mathrm{energy}}$ =  1.207E-05 & 
  $\Delta_{\mathrm{energy}}$ =  1.173E-05 & 
  $\Delta_{\mathrm{energy}}$ =  1.179E-05 & 
  $\Delta_{\mathrm{energy}}$ =  1.193E-05  
  \\
  $\#\Th = 2\,248$, $\DoF$ = 107\,904 &
  $\#\Th = 247$, $\DoF$ =11\,856  &
  $\#\Th = 247$, $\DoF$ =11\,856&
  $\#\Th = 319$,  $\DoF$ =15\,312 \\  
  \end{tabular}
  }
  \end{center}
  \caption{Smooth rising thermal bubble: meshes (first line), isolines of
    the potential temperature perturbation with the contour interval of $0.05$
    (second line) and the isolines of the density
    with the contour interval of $0.05$
    (third line); the computations using the uniform grid (left column)
    and anisotropic meshes generated with the tolerances
    $\TOL= 0.02$, $0.015$ and $0.01$ (from second to fourth columns) at $T=500\si{s}$.
    Note: The numbers of $\#\Th$ and $\DoF$ are the same for $\TOL=0.02$ and
   $\TOL=0.015$ which is just a coincidence.} 
  \label{fig:ANGENER}
\end{figure}

The goal of the following numerical experiments is to demonstrate that
(i) the use of different meshes at different time layers does not have a significant
impact on the accuracy of the computations and
(ii) the density is suitable as the quantity of interest in the {\hpAMA} algorithm.
We consider the rising thermal bubble test case
in Section~\ref{sec:thermal} with the initial condition \eqref{thermal1}.
We compare results achieved on
\begin{itemize}
\item quasi-uniform grid having
2\,248 elements using $P_3$ polynomial approximation with respect to space and
$P_1$ polynomial approximation with respect to time variables,
\item anisotropic $hp$-meshes  produced by Algorithm~\ref{alg:hp-AMA-ST}
  using density as the quantity $\qht$ with three tolerances $\TOL$.
\end{itemize}
Figure~\ref{fig:ANGENER} shows the final mesh (without the polynomial degree distribution),
the isolines
of the potential temperature perturbation  and the isolines of the density.
We observe that the potential temperature perturbation isolines give the
expected mushroom cloud with increasing accuracy for the increasing mesh refinement.
On the other hand, the density distribution is constant in the horizontal
direction and linear in the vertical one with a small perturbation
around the mushroom cloud. Even these small perturbations are sufficient for
the {\hpAMA} algorithm to detect the regions of the largest interpolation error
and using the evaluated higher order directional derivatives, the algorithm sets the
required shape and orientation of the mesh elements.

Moreover, we evaluate the quantities
\begin{align}
  \Delta_{\mathrm{mass}} = \left|\frac{M(T) - M(0)}{M(0)}\right|\ \ \mbox{with}\ \ 
  M(t) = \int_{\Om} \rho_{h\tau}(x,t)\dx
  \quad\mbox{and}\quad
  \Delta_{\mathrm{energy}} = \left|\frac{E(T) - E(0)}{E(0)}\right|\ \ \mbox{with}\ \ 
  E(t) = \int_{\Om} \ener_{h\tau}(x,t)\dx,
\end{align}
where $\rho_{h\tau}$ and $\ener_{h\tau}$ are the approximations of the density and the energy,
respectively, given by the STDG method.
The corresponding values are given at the bottom of Figure~\ref{fig:ANGENER} together
with the number of mesh elements $\Thm$ and $\DoF=\Nhm=\dim\bSmhp$ for the last time layer $m=r$.
We observe the same level of violation of the mass and energy conservation
for the computations using fixed and varying meshes
which indicates that the resulting numerical solutions do not suffer from
the errors arising from multiple recomputations between
adapted meshes.

\section{Numerical experiments}
\label{sec:num}

In this section, we present the numerical solution of several benchmarks
from \cite{GiraldoRostelli_JCP08} by the combination of the 
STDG and {\hpAMA} methods, cf.~Algorithm~\ref{alg:hp-AMA-ST}.
The computations have been carried out by the in-house code ADGFEM \cite{ADGFEM}
employing ANGENER code \cite{angener} for the mesh adaptation.
For some of the presented examples,
the {\hpAMA} method significantly reduces the computational costs
in comparison to non-refined grids. On the other hand, some test problems do not have anisotropic
features so the {\hpAMA} method generates grids that are close to the uniform ones.
Similarly as in \cite{GiraldoRostelli_JCP08}, all examples except
the sharp thermal bubble in Section~\ref{sec:thermal} and 
the density current flow in
Section~\ref{sec:density} are treated as inviscid flows, i.e., $\bRv=0$ in \eqref{eq:NS}.
The physical constants appearing in \eqref{eq:w} -- \eqref{eq:tempP} take the values
$\kappa=1.4$, $\Pr=0.72$ and $\ccV=718\si{J\cdot kg^{-1} \cdot K^{-1}}$.
We use the piecewise linear approximation with respect to time ($q=1$ in \eqref{Shpq}) 
for all numerical examples except the steady-state problem in Section~\ref{sec:schar} when
we set $q=0$.

\subsection{Inertia gravity wave}

The inertia gravity wave benchmark
(cf.~\cite{GiraldoRostelli_JCP08,SkamarockKlemp_MWR94})
exhibits the evolution of a potential temperature perturbation
in the rectangular domain $\Om=(0,\np{300000})\times (0,\np{10000})\si{m}$
with the periodic boundary conditions on the left and right boundaries of $\Om$
and the no-flux boundary conditions on the horizontal boundary parts.
The initial state of the atmosphere has the mean horizontal velocity $v_1=20\si{m\cdot s^{-1}}$
in a uniformly stratified atmosphere with the Brunt-V\"ais\"al\"a
frequency $N=0.01\si{s^{-1}}$ where the mean potential temperature  fulfills
\begin{align}
  \label{IG1}
  \bar{\tempP} = \tempP_0 \exp(N^2 x_2 / g),\qquad \tempP_0 = 300 \si{K}.
\end{align}
The Exner pressure is then given by the hydrostatic balance
equation \eqref{balance}, namely
\begin{align}
  \label{IG2}
  \bar{\pressE} =1 + \frac{g^2}{\ccP\,\tempP_0\, N^2}\left(\exp(-N^2 x_2 / g) -1\right).
\end{align}
Moreover, to the mean potential temperature $\bar{\tempP}$, we add the perturbation
\begin{align}
  \label{IG3}
  {\tempP}' = \tempP_c \sin(\pi x_2 / h_c) \left(1 + ((x_1-\bar{x}_1)/a_c)^2\right)^{-1},
\end{align}
where $\tempP_c= 0.01\si{K}$, $h_c = \np{10000}\si{m}$, $a_c = \np{5000}\si{m}$
and $\bar{x}_1 = \np{100000}\si{m}$.
The initial perturbation propagates to the left and right symmetrically,
but it is carried away by the mean horizontal flow till
the final time $T=3\,000\si{s}$.

\begin{figure}
  \includegraphics[width=0.490\textwidth]{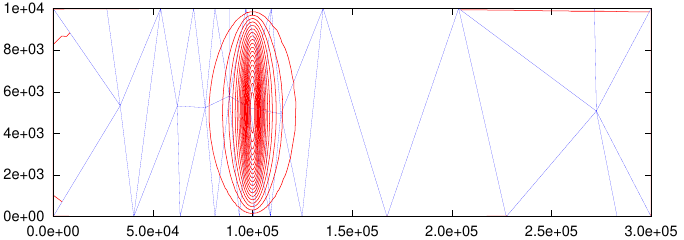}
  \hspace{0.01\textwidth}
  \includegraphics[width=0.490\textwidth]{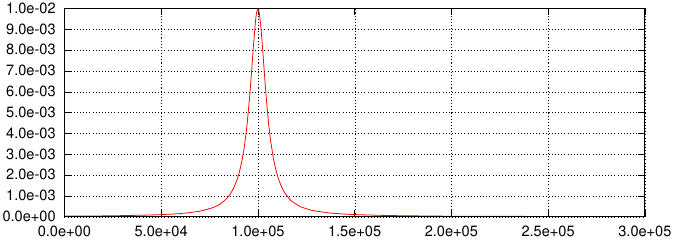}
  \\[3mm]
  \includegraphics[width=0.490\textwidth]{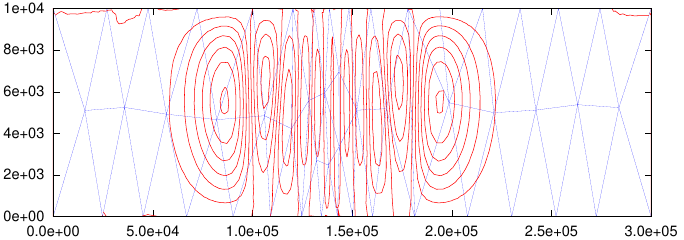}
  \hspace{0.01\textwidth}
  \includegraphics[width=0.490\textwidth]{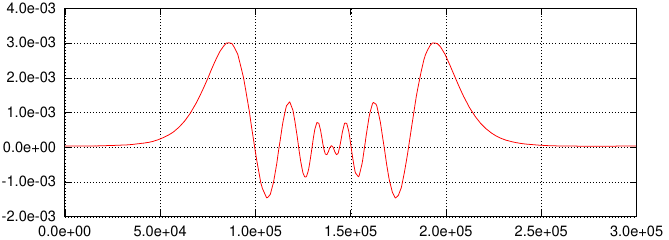}
  \\[3mm]
  \includegraphics[width=0.490\textwidth]{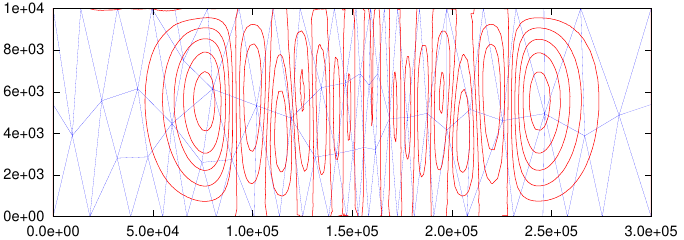}
  \hspace{0.01\textwidth}
  \includegraphics[width=0.490\textwidth]{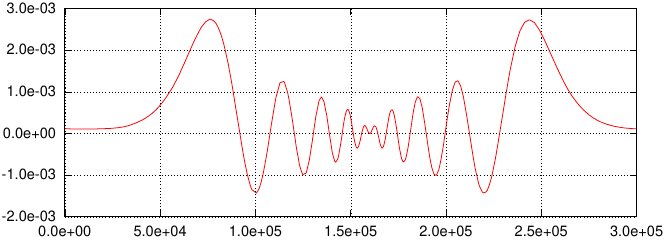}
  \caption{Interior gravity wave: isolines of the potential temperature perturbation
    (the contour interval of $5\cdot10^{-4}\si{K}$) with the corresponding mesh
    (left) and the potential temperature perturbation along
    $x_2=5\, 000\si{m}$
    at $t=0\si{s}$ (top), $t=\np{2000}\si{s}$ (center) and $t=\np{3000}\si{s}$ (bottom).}
  \label{fig:IGW1}
\end{figure}

The solution of this problem is smooth, so we expect that the mesh adaptive
algorithm generates a sequence of coarse meshes with a low number of elements
and with high polynomial approximation degrees.
Figure~\ref{fig:IGW1} presents the potential temperature perturbation at the physical
times $t=\np{0}\si{s}$, $t=\np{2000}\si{s}$ and $t=\np{3000}\si{s}$,
namely the isolines and the horizontal distribution
along $x_2=\np{5000}\si{m}$. The left figures indicate also the generated triangular meshes,
the polynomial approximation degrees $\pK$ for (almost) all $K\in\Thm$
are equal to 8 which is the maximal degree supported by our code.

We observe a smooth resolution of the isolines although relatively coarse meshes have been
generated. The profiles are perfectly symmetric along the vertical position
$x_1=\bar{x}_1+ v_1 t$, i.e., $x_1= \np{140000}\si{m}$ for $t=\np{2000}\si{s}$
and $x_1= \np{160000}\si{m}$ for $t=\np{3000}\si{s}$.
The results are in full agreement with the reference data in
\cite[Section~6.1]{GiraldoRostelli_JCP08}.

\subsection{Rising thermal bubbles}
\label{sec:thermal}

The rising thermal bubble test case simulates the evolution of a warm bubble
in a constant potential temperature field.
Since the bubble is warmer,  it rises while deforming as a consequence of
the shearing motion caused by the velocity field gradients until it forms a mushroom cloud.
According to \cite{GiraldoRostelli_JCP08,Robert_JAS93},
we consider the computational domain $\Om=(0,1000)\times(0,1000)\si{m}$,
the constant mean potential temperature $\bar{\tempP} = 300\si{K}$ and
the Exner pressure follows from the hydrostatic balance \eqref{balance}
as $\pressE = 1 - g/(\ccP \bar{\tempP}) x_2$. The initial velocity is equal to zero.
Moreover, to the mean flow, we add two different  potential temperature perturbations,
the smooth one
\begin{align}
  \label{thermal1}
  \tempP' = 
  \begin{cases}
    0 & \mbox{for } r > r_c, \\
    \tfrac12 \tempP_c \left(1 + \cos(\pi\, r / r_c)\right)    & \mbox{for } r \le r_c, \\
  \end{cases}
\end{align}
where $\tempP_c= 0.5\si{K}$, $r = | x - x_c|$, $x_c=(500, 300)\si{m}$ and $r_c = 250\si{m}$,
and the sharp one
\begin{align}
  \label{thermal2}
  \tempP' = 
  \begin{cases}
    \tempP_c & \mbox{for } r < r_c, \\
    \tempP_c  \exp( -(r-r_c)^2/s^2  ) & \mbox{for } r \ge r_c, \\
  \end{cases}
\end{align}
where $\tempP_c= 0.5\si{K}$, $r = | x - x_c|$, $x_c=(500, 260)\si{m}$, $r_c = 150\si{m}$ and
$s=10$.
We note that in the original paper \cite{Robert_JAS93},
the parameter $s$ in \eqref{thermal2} was chosen as $s=100$, so we have much sharper initial
perturbation in our case.
For the smooth initial bubble, the final physical time is set to $T=700\si{s}$.
For the sharper one, we consider $T=850\si{s}$ and a small amount of the dynamic viscosity
$\mu=3.44\cdot10^{-3} \si{Pa\cdot s}$. 
The no-flux boundary conditions are prescribed on the whole boundary.

\begin{figure}
  \begin{center}

  \includegraphics[width=0.325\textwidth]{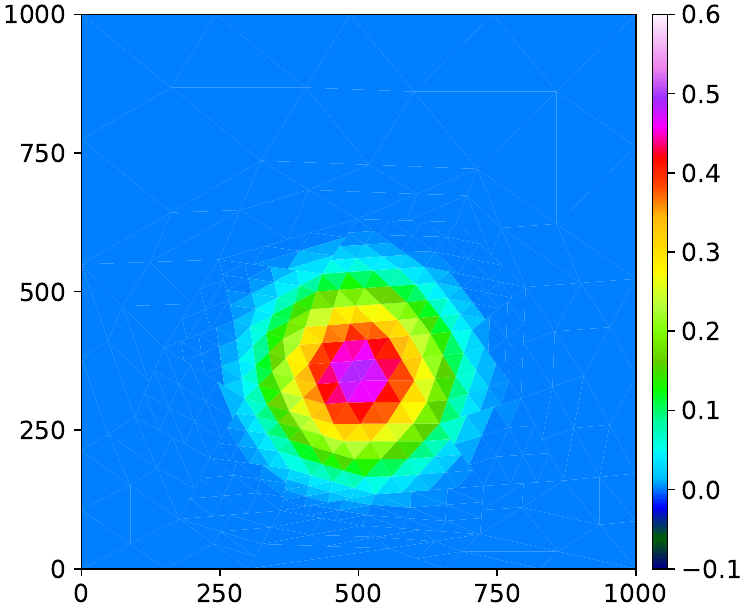}
  \hspace{0.000\textwidth}
  \includegraphics[width=0.325\textwidth]{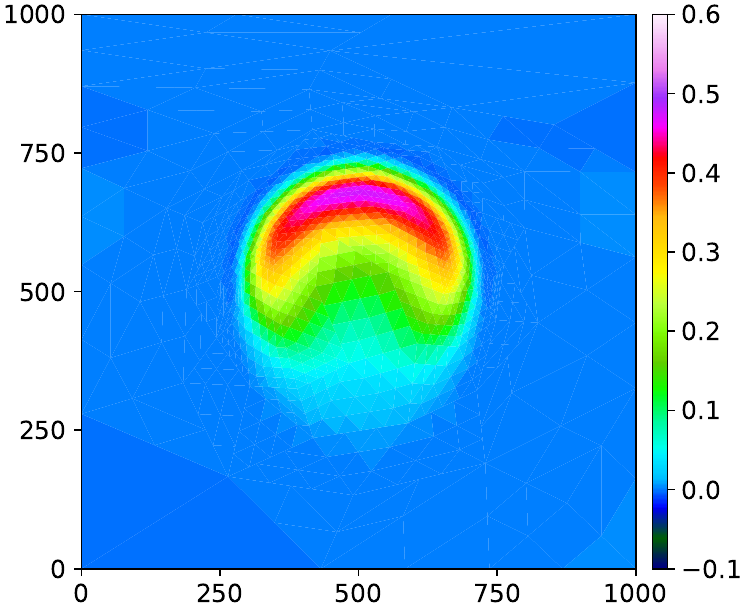}
  \hspace{0.000\textwidth}
  \includegraphics[width=0.325\textwidth]{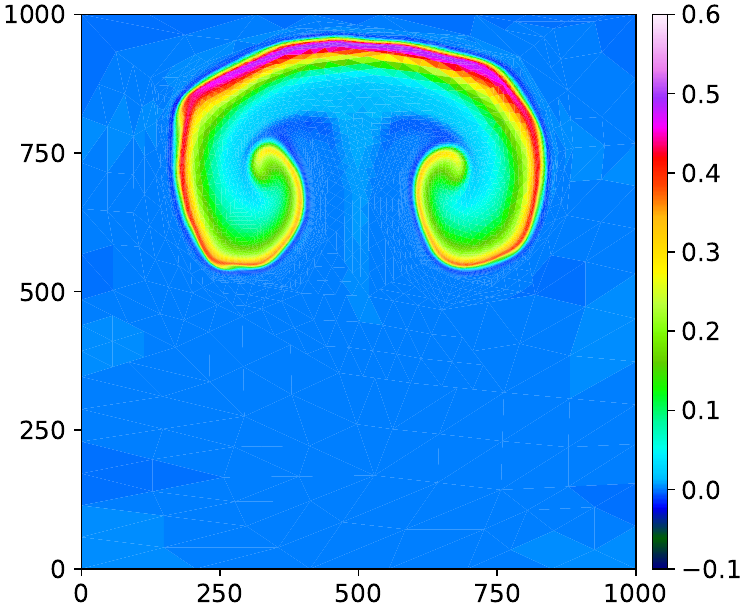} \\[2mm]
  
  \includegraphics[width=0.325\textwidth]{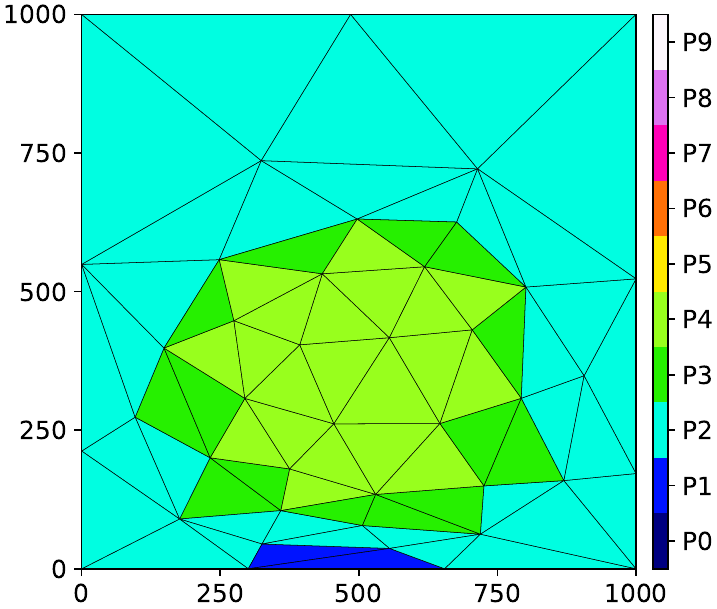}
  \hspace{0.000\textwidth}
  \includegraphics[width=0.325\textwidth]{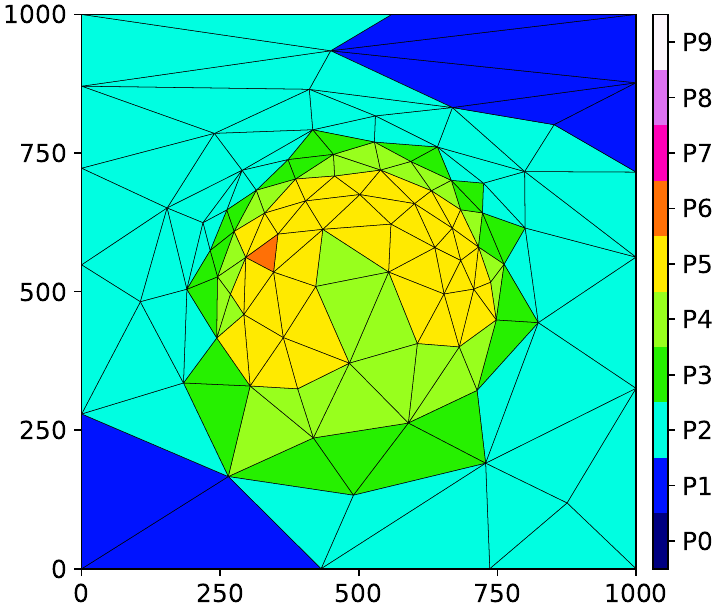}
  \hspace{0.000\textwidth}
  \includegraphics[width=0.325\textwidth]{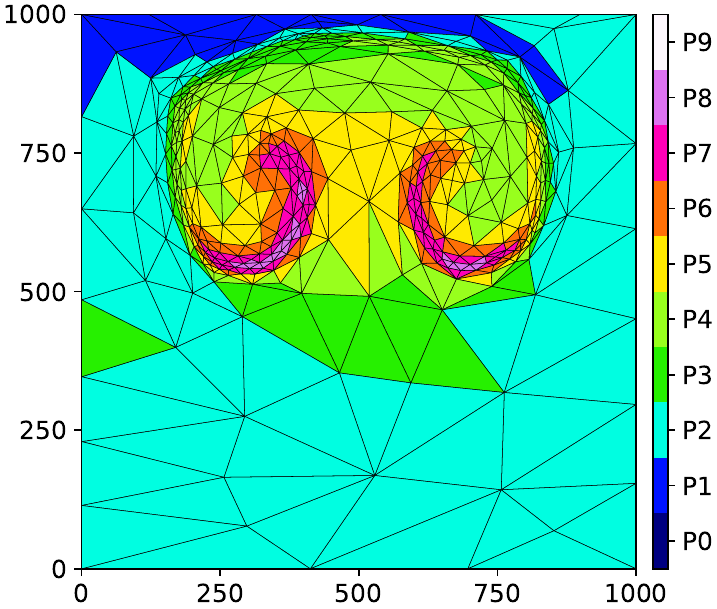}
  
  \end{center}
  \caption{Smooth rising thermal bubble: the distribution of the potential temperature perturbation
    (top) and the corresponding $hp$-meshes (bottom)
    at $t=0\si{s}$, $t=350\si{s}$ and $t=700\si{s}$ (from left to right).}
  \label{fig:thermalS}
\end{figure}

\begin{figure}
  \begin{center}

  \includegraphics[width=0.325\textwidth]{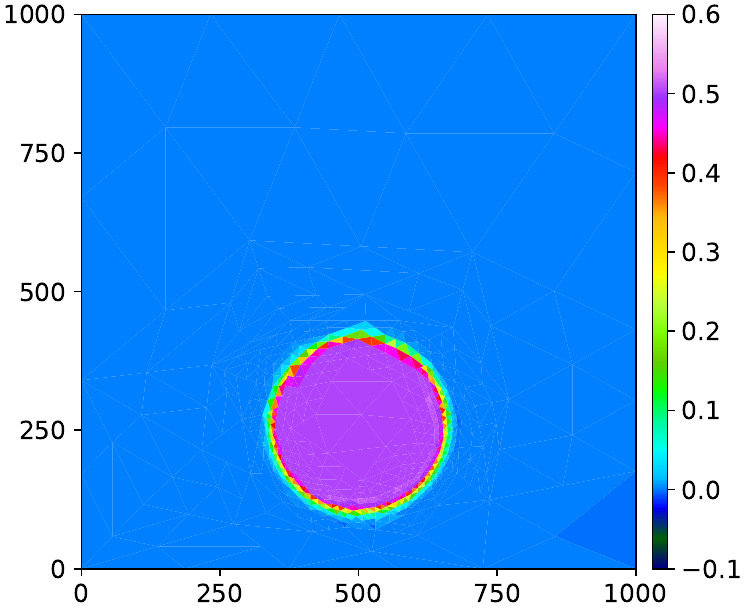}
  \hspace{0.000\textwidth}
  \includegraphics[width=0.325\textwidth]{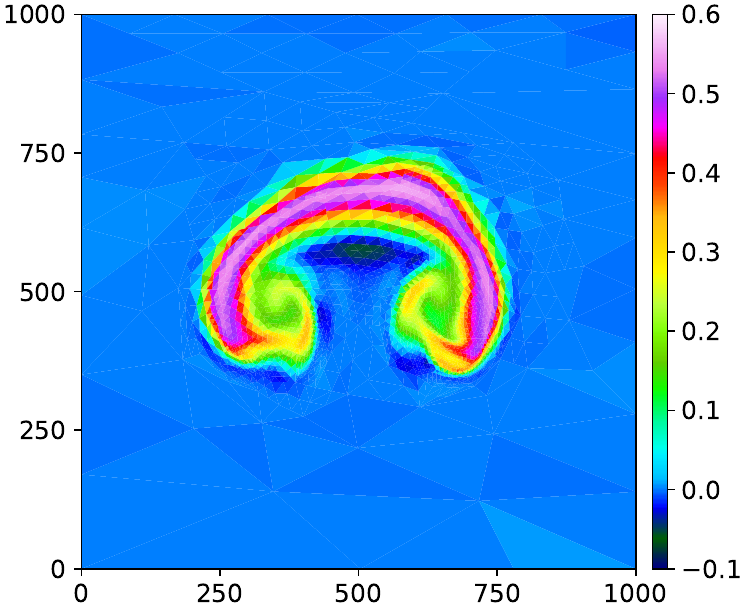}
  \hspace{0.000\textwidth}
  \includegraphics[width=0.325\textwidth]{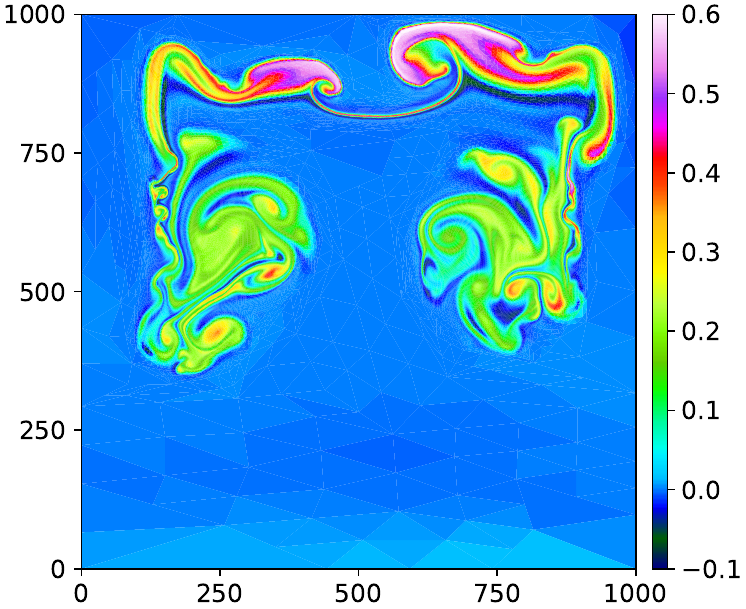} \\[2mm]
  
  \includegraphics[width=0.325\textwidth]{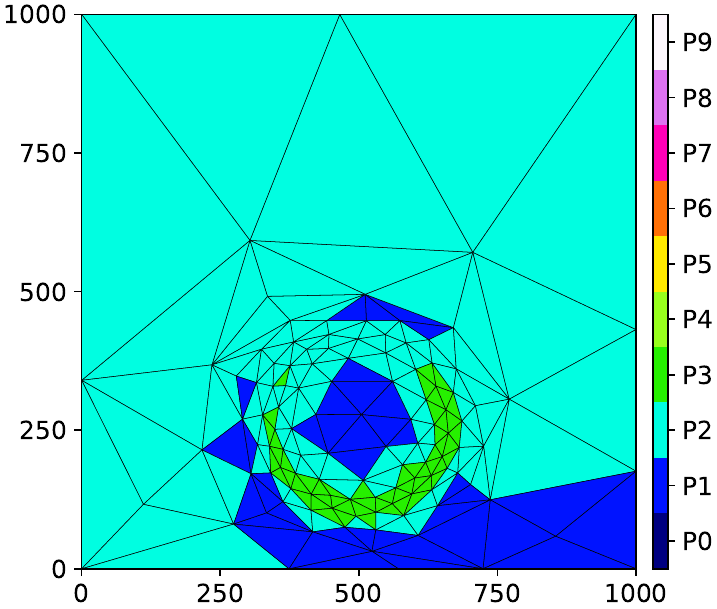}
  \hspace{0.000\textwidth}
  \includegraphics[width=0.325\textwidth]{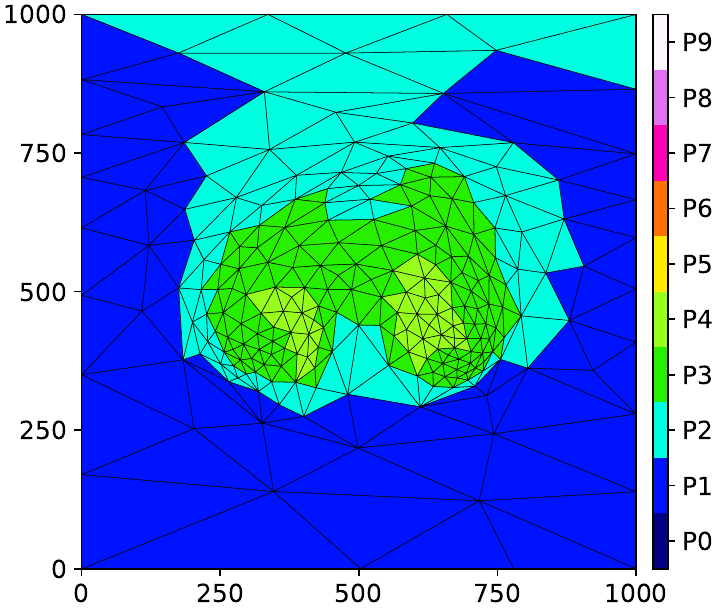}
  \hspace{0.000\textwidth}
  \includegraphics[width=0.325\textwidth]{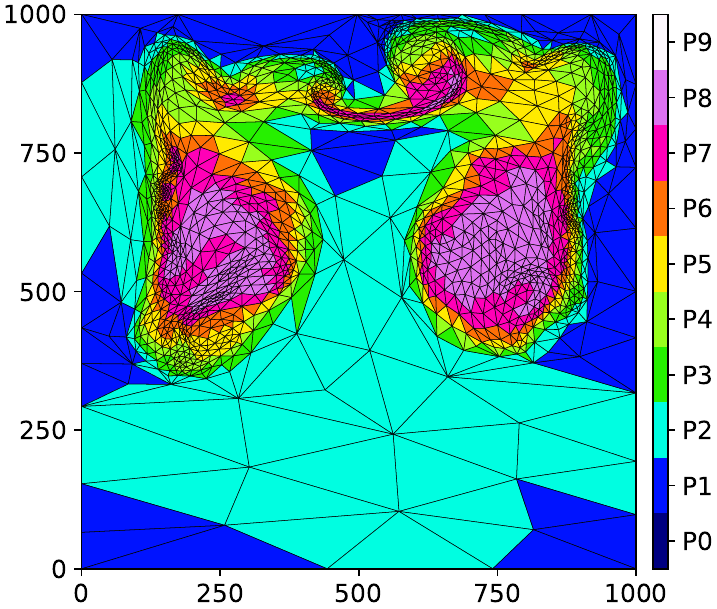}
  
  \end{center}
  \caption{Sharp rising thermal bubble: the distribution of the potential temperature perturbation
    (top) and the corresponding $hp$-meshes (bottom)
    at $t=0\si{s}$, $t=425\si{s}$ and $t=850\si{s}$ (from left to right).}
  \label{fig:thermalF}
\end{figure}

Figures~\ref{fig:thermalS} and \ref{fig:thermalF} show the
distribution of the potential temperature perturbation
and the corresponding $hp$-meshes at selected time steps for both initial settings.
We observe an accurate capture of the deformed bubbles together with the expected
anisotropic mesh refinement. We note that all figures are visualized by
the library Matplotlib \cite{Matplotlib} where no additional smoothing of the results
is done.

\subsection{Density current flow problem}
\label{sec:density}

The density current flow problem benchmark  was proposed in \cite{StrakaALL_IJNMF93}
and it exhibits the evolution of a cold bubble dropped in a neutrally stratified atmosphere.
Since the bubble is cold, it sinks to the ground and starts to move along the ground
forming Kelvin--Helmholtz vortexes.

In the same way as \cite{GiraldoRostelli_JCP08},
we consider the computational domain $\Om=(0,\np{25600})\times(0,\np{6400})\si{m}$, the final time
$T=900\si{s}$. We prescribe the symmetric boundary condition along $x_1=0$ and the
no-flux boundary conditions over the rest of the boundary.
It was discussed in \cite{StrakaALL_IJNMF93} that a viscosity is required to
obtain a grid-converged solution, hence we set $\mu = 0.1\si{Pa\cdot s}$ in \eqref{eq:tau}.
The initial flow is given by the velocity $\bkv=0$
and the constant mean potential temperature
$\bar{\tempP} = 300\si{K}$  perturbed  by
\begin{align}
  \tempP' = \frac{\tempP_C}{2}\left( 1 + \cos(\pi r)\right),\qquad \tempP_C= - 15\si{K},
\quad  r = \left( \left({(x_1 - x_1^c)}/{x_1^r}\right)^2
  + \left( {(x_2 - x_2^c)}/{x_2^r}\right)^2
  \right)^{1/2},
\end{align}
where the center of the bubble is at $(x_1^c, x_x^c) = (0, \np{3000})\si{m}$ and
its size is $(x_1^r,x_2^r) = \left(\np{4000},\np{2000}\right)\si{m}$.
The Exner pressure follows from the balance equation \eqref{balance}.

\begin{figure}
  \begin{center}

  \includegraphics[width=0.495\textwidth]{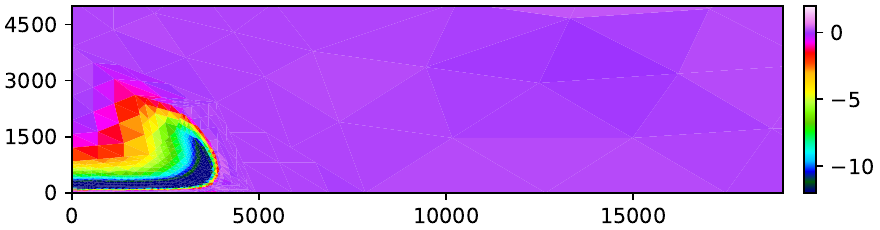}
  \hspace{0.000\textwidth}
  \includegraphics[width=0.495\textwidth]{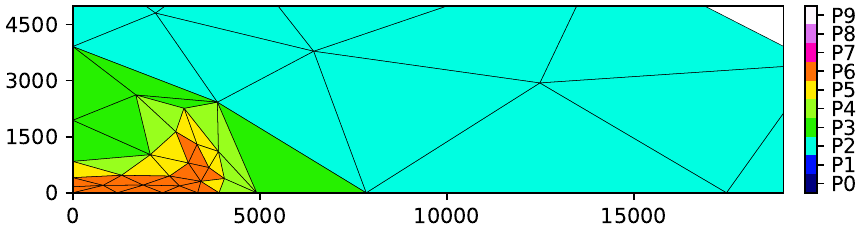}

  \vspace{2mm}

  \includegraphics[width=0.495\textwidth]{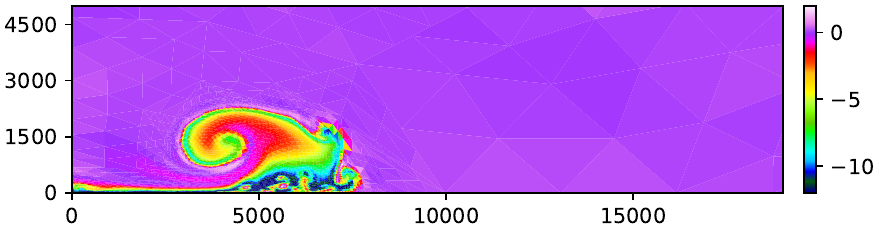}
  \hspace{0.000\textwidth}
  \includegraphics[width=0.495\textwidth]{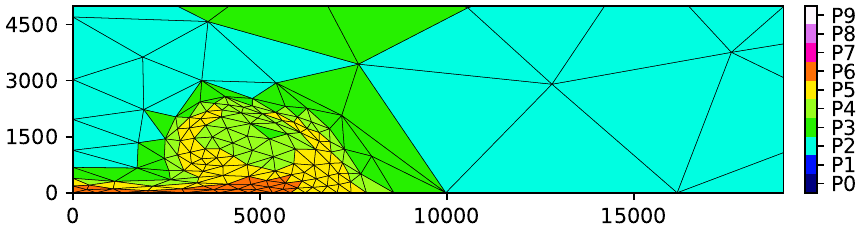}

  \vspace{2mm}

  \includegraphics[width=0.495\textwidth]{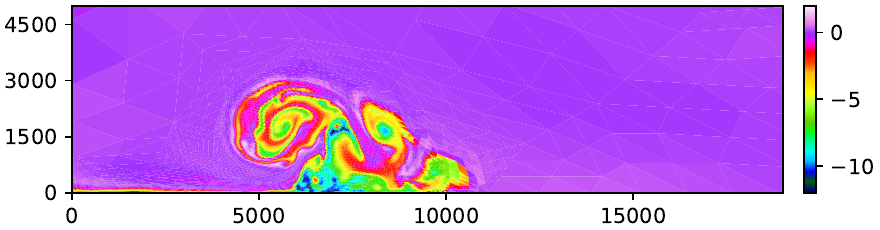}
  \hspace{0.000\textwidth}
  \includegraphics[width=0.495\textwidth]{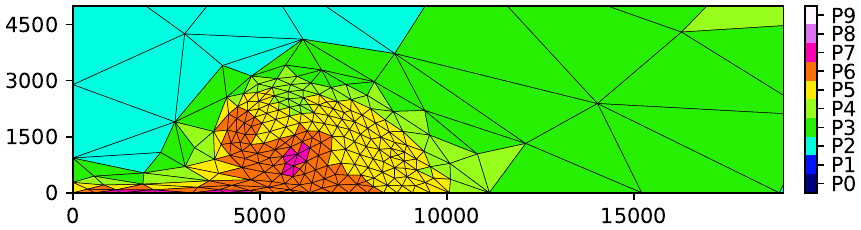}

  \vspace{2mm}
  
  \includegraphics[width=0.495\textwidth]{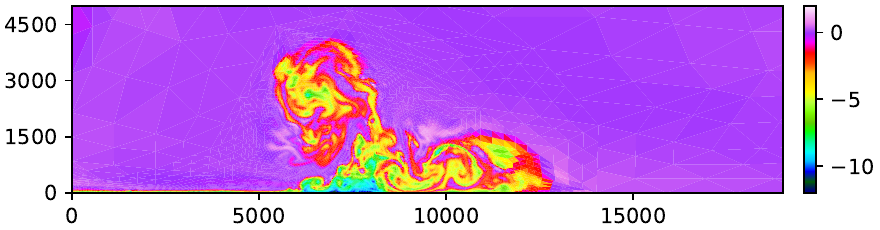}
  \hspace{0.000\textwidth}
  \includegraphics[width=0.495\textwidth]{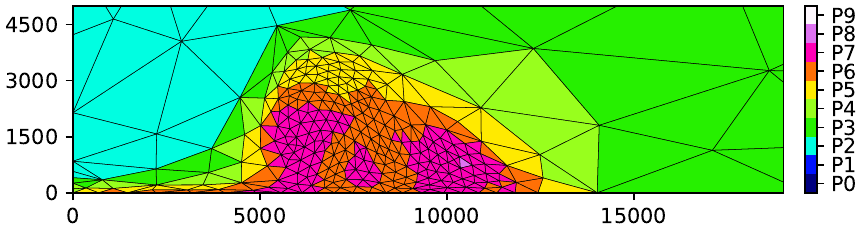}
  
  \end{center}
  \caption{Density current flow: the distribution of the potential temperature perturbation
    (left) and the corresponding $hp$-meshes (right)
    at $t=300\si{s}$, $t=500\si{s}$ , $t=700\si{s}$ and $t=\np{900}\si{s}$ (from top to bottom).}
  \label{fig:density}
\end{figure}

Figure~\ref{fig:density} shows the distribution of the perturbation of the
potential temperature and the corresponding $hp$-meshes at selected time levels.
The development of Kelvin--Helmholtz vortexes is easily observed together with
the expected mesh adaptation. Due to high polynomial approximation degrees, 
the vortexes are captured sharply even for relatively coarse meshes, the maximal resolution
is approximately $250\si{m}$. We note that in the reference paper \cite{StrakaALL_IJNMF93},
the larger viscosity $\mu = 75\si{Pa\cdot s}$ has been used.

\subsection{Sch\"ar mountain}
\label{sec:schar}

This test case deals with the steady-state solution of hydrostatic flow over a five-peak mountain
chain, with steady inflow and outflow boundary conditions.
In the same way as in \cite{GiraldoRostelli_JCP08,ScharALL_MWR02},
the initial state of the atmosphere has a constant mean flow $v_1= 10\si{m\cdot s^{-1}}$
in a uniformly stratified atmosphere with a Brunt-V\"ais\"al\"a frequency of $N=0.01\si{s^{-1}}$.
Using \eqref{IG1}--\eqref{IG2}, we set the reference potential
temperature and the Exner pressure with $\tempP_0 = 280 \si{K}$.
The computational domain is defined by
\begin{align}
  \label{sch1}
  \Om&=\left\{(x_1,x_2),\   g(x_1) < x_2 < \np{21000},\ \np{-25000} < x_1 < \np{25000}
  \right\}, \quad
   g(x_1)= h_c\exp(-(x_1/a_c)^2)\cos^2(\pi\,x_1/\lambda_c), 
\end{align}
where 
with the parameters $h_c=\np{250}\si{m}$, $\lambda_c= \np{4000}\si{m}$ and $a_c = \np{5000}\si{m}$.
No-flux boundary conditions are used along the bottom surface
$\{(x_1, g(x_1)),  x_1\in(\np{-25000}, \np{25000})\}$,
while
non-reflecting boundary conditions are used on the rest of boundary $\gom$.

Since the flow is steady, we adopt Algorithm~\ref{alg:hp-AMA-ST}.
In step \ref{alg:G1w},
we set $\etam:=\norm{\Wqht(t_m) -\Pi\Wqht(t_m)}{L^2(\Om)}$ and
carry out the computation formally for $T\to \infty$.
Moreover, we do not repeat the time steps
if the stopping criterion (step \ref{alg:G1z})  is not met and proceed directly to the next time
step. Furthermore, the constant $c_T$ in \eqref{CT} is set large since we do not take 
care of the time discretization error, hence we replace \eqref{CA} by
$\etaAm(\wht) \leq c_A \etaSm(\wht)$. The computation is stopped when 
$\etam\leq\TOL$.

\begin{figure}
  \includegraphics[width=0.490\textwidth]{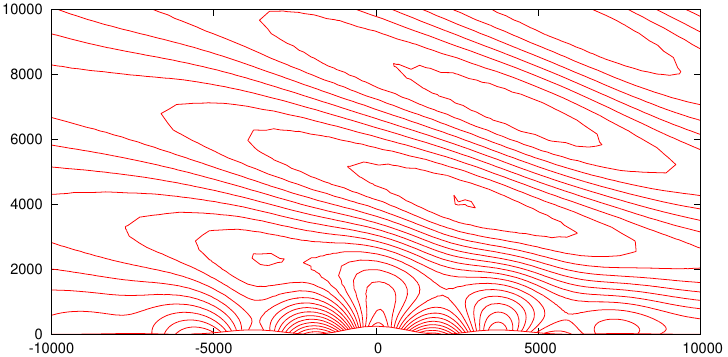}
  \hspace{0.01\textwidth}
  \includegraphics[width=0.490\textwidth]{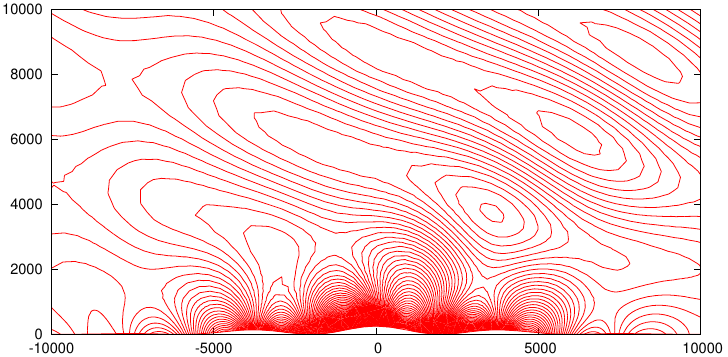}
  \caption{Sch\"ar mountain: zoom of isolines of the horizontal (left) and vertical (right)
    components of velocity 
    with the contour interval of $0.2$ and $0.0005$, respectively.}
  \label{fig:Schar1}
\end{figure}

\begin{figure}
  \includegraphics[height=0.2\textwidth]{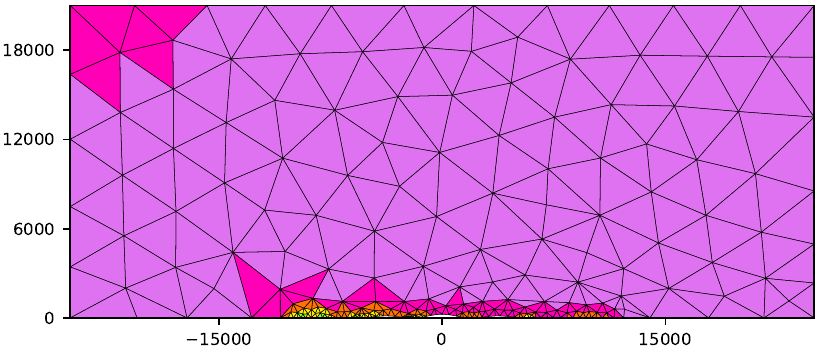}
  \hspace{0.01\textwidth}
  \includegraphics[height=0.2\textwidth]{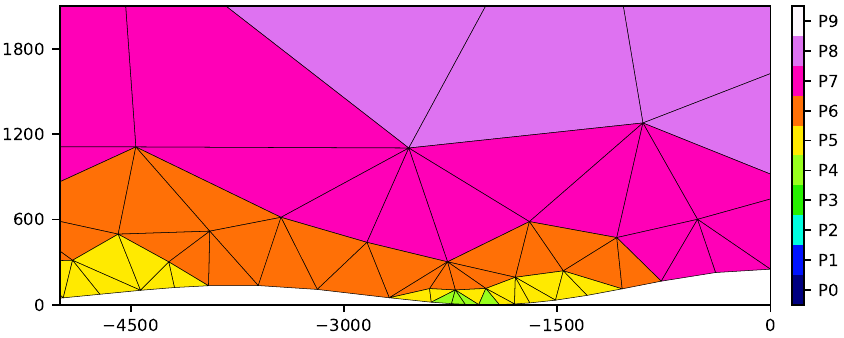}
  \caption{Sch\"ar mountain: the final $hp$-mesh, total view (left) and
    its zoom around the bottom part of the boundary (right).}
  \label{fig:Schar2}
\end{figure}

Figure~\ref{fig:Schar1} shows  isolines of the horizontal  and vertical
components of the velocity after 20 levels of the mesh adaptation when the tolerance
$\TOL=10^{-5}$ has been achieved. 
The corresponding $hp$-mesh with a zoom of the bottom boundary part
is presented in Figure~\ref{fig:Schar2}.
We note that
the mesh triangles adjacent to the non-polygonal part of the boundary are treated
as the curvilinear elements, namely we use a piecewise cubic approximation.
We observe small elements near the bottom part of the boundary which are 
necessary to capture its geometry. On the other hand,
the mesh is mostly quasi-uniform with large elements having the highest
polynomial approximation degree
in the majority of the domain since this example does not contain any strong anisotropic
feature. Again, these results are in agreement with the reference data in
\cite[Section~6.4]{GiraldoRostelli_JCP08}.


\section{Conclusion}
\label{sec:concl}

We presented the $hp\tau$-adaptive space-time discontinuous Galerkin method for the 
numerical solution of non-hydrostatic mesoscale atmospheric problems.
The time discontinuous discretization with a weak coupling by the penalty
admits using different
meshes at different time levels in a natural way. We demonstrated by numerical experiments
that the complete re-meshing during the computational process does not significantly affect
the accuracy of the approximation. Moreover, measuring the interpolation error
in terms of the density is sufficient for the anisotropic mesh adaptation algorithm 
to generate appropriately adapted meshes, namely to capture physically relevant quantity as
the perturbation of the potential temperature. The robustness of the adaptive algorithm
is supported by numerical simulation of several benchmarks.



\end{document}